\documentclass[10pt,twoside]{siamart1116}

\usepackage[english]{babel}
\usepackage{graphicx,epstopdf,epsfig}
\usepackage{amsfonts,epsfig,fancyhdr,graphics, hyperref,amsmath,amssymb}

\usepackage{mathdots}
\usepackage{algorithm}
\usepackage{algpseudocode}

\newcommand{\HE}{Name of Handling Editor}
\newcommand{\DoS}{Month/Day/Year}
\newcommand{\DoA}{Month/Day/Year}
\newcommand{\CA}{Heike Fa\ss bender}

\newcommand{\Names}{Christian Bertram, Heike Fa\ss bender}
\newcommand{\Title}{On a family of low-rank algorithms for large-scale algebraic Riccati equations}

\newtheorem{remark}[theorem]{Remark}

\renewcommand{\theequation}{\thesection.\arabic{equation}}
\renewtheorem{theorem}{Theorem}[section]

\numberwithin{equation}{section} 

\begin{document}

\bibliographystyle{plain}

\setcounter{page}{1}

\thispagestyle{empty}

 \title{\Title\thanks{Received
 by the editors on \DoS.
 Accepted for publication on \DoA. 
 Handling Editor: \HE. Corresponding Author: \CA}}

\author{
Christian Bertram\thanks{
Institute for Numerical Analysis, TU Braunschweig, Universit\"atsplatz 2, 38106 Braunschweig, Germany}
\and 
Heike Fa{\ss}bender\thanks{
Institute for Numerical Analysis, TU Braunschweig, Universit\"atsplatz 2, 38106 Braunschweig, Germany (h.fassbender@tu-braunschweig.de)}}

\markboth{\Names}{\Title}

\maketitle

\begin{abstract}
In \cite{BenBKS18} it was shown that four seemingly different algorithms for computing low-rank approximate solutions $X_j$ to the solution $X$ of large-scale continuous-time algebraic Riccati equations (CAREs) $0 = \mathcal{R}(X) := A^HX+XA+C^HC-XBB^HX $ generate the same sequence $X_j$ when used with the same parameters.
The  Hermitian low-rank approximations $X_j$ are of the form $X_j = Z_jY_jZ_j^H,$ where $Z_j$ is a matrix with only few columns and $Y_j$ is a small square Hermitian matrix. Each $X_j$ generates a low-rank Riccati residual $\mathcal{R}(X_j)$  such that the norm of the residual can be evaluated easily allowing for an efficient termination criterion.
Here a new family of methods to generate such low-rank approximate solutions $X_j$ of CAREs is proposed.
Each member of this family of algorithms proposed here generates the same sequence of $X_j$ as the four previously known algorithms. The approach is based on a block rational Arnoldi decomposition and an associated block rational Krylov subspace spanned by $A^H$ and $C^H.$   Two specific versions of the general algorithm will be considered; one will turn out to be a rediscovery of the RADI algorithm, the other one allows for a slightly more efficient implementation compared to the RADI algorithm (in case the Sherman-Morrision-Woodbury formula and a direct solver is used to solve the linear systems that occur). Moreover, our approach allows for adding more than one shift at a time.
\end{abstract}



\section{Introduction}
Finding the unique stabilizing solution $X = X^H$ of large-scale algebraic Riccati equations
\begin{equation}\label{ric0}
0 = \mathcal{R}(X) := A^HX+XA+C^HC-XBB^HX 
\end{equation}
with a large, sparse matrix $A \in \mathbb{C}^{n \times n},$ and matrices $B\in \mathbb{C}^{n \times m}$ and $C\in \mathbb{C}^{p \times n}$ is of interest in a number of applications as noted in \cite{BenS13,LinS15,SimSM14} and references therein. Here, $B$ and $C$ are assumed to have full column and row rank, resp.,  with $m,p \ll n.$ Further, we assume that the unique stabilizing solution $X = X^H,$ which is positive semidefinite and makes the closed-loop matrix $A-BB^HX$ stable, exists. It exists if rank$[A-\lambda I ~~~ BB^H]=$ rank$[A^H-\lambda I ~~~ CC^H]=n$ for all $\lambda$ in the closed right half plane \cite{LanR95}. 
Even though $A$ is large and sparse, the solution $X$ will still be a dense matrix in general. But our assumptions on $B$ and $C$ often imply that the sought-after solution $X$ will have a low numerical rank (that is, its numerical rank is much smaller than $n$) \cite{BenB16}. This allows for the construction of iterative methods that approximate $X$ with a series of low-rank matrices stored in low-rank factored form. 

To be precise, we are interested in Hermitian low-rank approximations $X_j$ to $X$  of the form  $X_j = Z_jY_j Z_j^H,$ where $Z_j \in \mathbb{C}^{n \times {k_j}}$ is a rectangular matrix with only few columns ($k_j \ll n$) and $Y_j \in \mathbb{C}^{{k_j} \times {k_j}}$ is a small square Hermitian matrix. 
Any method which generates such low-rank approximations is especially suitable for large-scale applications as there is no need to store $X_j$ as a full dense matrix, but just the much smaller matrices $Z_j$ and $Y_j.$
There are several methods (e.g., rational Krylov subspace methods, low-rank Newton-Kleinman methods and Newton-ADI-type methods)
which produce such a low-rank approximation; see, e.g. \cite{AmoB10,BenBKS18,BenHSW16,HeyJ09,LinS15,Pal19,Sim16,SimSM14,WonB05,WonB07} and \cite{BenBKS20} for an overview. Basically, all these methods use certain (rational) Krylov subspaces as approximation spaces (that is, the space spanned by the columns of $Z_j$).

Our contribution in this paper builds up on \cite{BenBKS18}. 
Whether any Hermitian matrix $X_0 \in \mathbb{C}^{n \times n}$ is a good approximation to the desired solution $X$ of \eqref{ric0} is usually measured via the norm of the Riccati residual $\mathcal{R}(X_0).$ The idea  pursued in \cite{BenBKS18} is 
given an approximation $X_0$, determine an approximation $\widetilde E$ to $E$ with $\mathcal{R}(X_0+E)= 0,$ set $X_1 = X_0+\widetilde E,$ and repeat the process with $X_1$.
The authors show that, starting the proposed iteration with $X_0=0$ or any other $X_0$ such that $\mathcal{R}(X_0) = R_0R_0^H$ with a full rank matrix $R_0 \in \mathbb{C}^{n \times p},$ all subsequent updated approximations $X_j$  yield a Riccati residual $\mathcal{R}(X_j) = R_jR_j^H$ with a matrix $R_j \in \mathbb{C}^{n \times p}$ (see Proposition 1 and (12) in \cite{BenBKS18}). In other words, the Riccati residual is always of rank at most $p.$ This allows for an easy update of $\| \mathcal{R}(X_j)\|_F$ as instead of the $n \times n$ matrix  $R_jR_j^H$ only a $p \times p$ matrix $R_j^HR_j$ has to be considered, $\| R_jR_j^H\|_F = \|R_j^HR_j\|_F$. Moreover, the factor $R_j$ can be computed efficiently by an additive update from $R_{j-1}.$ 
 Furthermore, the approximate solutions $X_j=Z_jY_jZ_j^H$ are of rank $jp.$ The factor $Z_j$ is constructed via an incremental update from $Z_{j-1}$ as $Z_j =[Z_{j-1} ~~~\hat Z_j]$ with $\hat Z_{j} \in \mathbb{C}^{n \times p},$ while for the matrix $Y_{j}$  an update of the form $\left[ \begin{smallmatrix}Y_{j-1}&\\ & \hat Y_j\end{smallmatrix}\right]$ with $\hat Y_j \in \mathbb{C}^{p \times p}$ holds. Hence, the resulting RADI method allows for an efficient way to store the approximations $X_j,$ as well as for an efficient way to calculate the factors of the approximations and the corresponding residuals cheaply (only one scalar product in the SISO case ($m=p=1$)).

 In addition, the $n \times p$ block columns of $Z_j = [\hat Z_1~~\hat Z_2~~\cdots~\hat Z_j]$ belong to the rational Krylov subspace
\begin{align}
\begin{split}\label{eq_L}
\mathcal{L}_j(A^H,C^H,\Sigma_j) &= \operatorname{blockspan}\{(A^H+\sigma_1I)^{-1}C^H, 
(A^H+\sigma_2I)^{-1}C^H, \\
&\qquad\qquad\qquad\qquad \ldots, (A^H+\sigma_jI)^{-1}C^H\}.
\end{split}
\end{align}
Here, $\Sigma_j = \{\sigma_1, \ldots, \sigma_j\}$ denotes the set of $j$ shifts $\sigma_i$ in the open left half plane $\mathbb{C}_{-}$ such that $A^H+\sigma_i I$ is nonsingular for $i=1,\ldots,j$ (that is, $-\Sigma_j \cap \Lambda(A) = \emptyset$ for the spectrum $\Lambda(A)$ of $A$).
If the $jp$ column vectors in $\mathcal{L}_j(A^H,C^H,\Sigma_j)$ are linearly independent (this implies that the $\sigma_i, i=1,\ldots,j$ are pairwise distinct) and
the columns of $Z_j \in \mathbb{C}^{n \times jp}$ represent a basis for the rational Krylov subspace $\mathcal{L}_j$, then $X_j= Z_jY_jZ_j^H$ holds for some matrix $Y_j \in \mathbb{C}^{jp \times jp}$  \cite[Proposition 2]{BenBKS18}. This decomposition is not unique in the sense, that for any nonsingular matrix $T \in \mathbb{c}^{jp \times jp}$ we have
\begin{equation}\label{eq:unique}
X_j = Z_j T T^{-1} Y_j T^{-H} T^H Z_j^H = \breve Z_j \breve Y_j \breve Z_j^H
\end{equation}
with $\breve Y_j = T^{-1} Y_j T^{-H}.$ The columns of $\breve Z_j$ form a different basis of $\mathcal{L}_j$ than those of $Z_j.$ In other words, once the basis of $\mathcal{L}_j$ is fixed, the decomposition of $X_j= Z_jY_jZ_j^H$ is unique.

Finally, \cite[Theorem 2]{BenBKS18} states that the approximation $X_j^\text{cay}$
of the Riccati solution obtained by the Cayley transformed Hamiltonian subspace iteration \cite{LinS15} and the approximation $X_j^\text{qadi}$ obtained by the qADI iteration \cite{WonB05,WonB07} are equal to the approximation $X_j^\text{radi}$ obtained by the RADI method
\begin{equation}\label{eq:equiv1}
X_j^\text{radi} = X_j^\text{cay} = X_j^\text{qadi}
\end{equation}
(if the initial approximation in all algorithms is zero and the same shifts are used). Beyond that 
\begin{equation}\label{eq:equiv2}
X_j^\text{radi} = X_j^\text{cay} = X_j^\text{qadi} = X_j^\text{inv}
\end{equation}
 if rank $C=1$ and the shifts are chosen equal to the distinct eigenvalues of $\left[\begin{smallmatrix} A & BB^H\\ C^HC & -A^H \end {smallmatrix}\right]$, where $X_j^\text{inv}$ is the approximation obtained by the invariant subspace approach \cite{AmoB10}. Parts of these connections have already been described in \cite{BenB16,LinS15}. From here on we will use the term Riccati ADI methods (proposed in \cite{BenBKS18}) to refer to these four equivalent methods.

Intrigued by the fact that the four above mentioned Riccati ADI methods
 produce (at least theoretically) the same approximate solution,  our first aim is to characterize all Hermitian rank-$jp$-matrices $X_j$ which 
can be written in the form
\[
X_j = Z_jY_jZ_j^H
\]
with $Z_j \in \mathbb{C}^{n \times jp}$  such that the $n \times p$ block columns of $Z_j$ span the block rational Krylov subspace 
\[ 
\operatorname{blockspan}(Z_j) = \frak{q}_j(A)^{-1}\operatorname{blockspan}\{C^H,A^HC^H, \ldots, (A^H)^{j-1}C^H\} = \mathfrak{K}_j^\Box(A^H,C^H,\frak{q}_j)
\]
for a polynomial $\frak{q}_j$ of degree $j$ and some nonsingular matrix $Y_j \in \mathbb{C}^{jp \times jp}$
and which yield
\[
\operatorname{rank}(\mathcal{R}(X_j)) = p.
\]
We will see that by choosing a space $\mathfrak{K}_j^\Box(A^H,C^H,\frak{q}_j)$, $X_j$ is  unique (just its low rank decomposition is only essentially unique in the sense of \eqref{eq:unique}). As
\begin{equation}\label{eq:same}
\mathcal{L}_j(A^H,C^H,\Sigma_j) = \mathfrak{K}_j^\Box(A^H,C^H,\frak{q})
\end{equation}
 in case the shifts $\sigma_i$ in $\Sigma_j$ are pairwise distinct and the negative roots of $\mathfrak{q}_j$, our results not only yield an elegant proof of \eqref{eq:equiv1} and \eqref{eq:equiv2}, but also enables the development of further, new methods generating the same sequence of approximations. Any such method will allow for an efficient way to store the approximations $X_j,$ as well as for an efficient way to calculate the quite reliable termination criterion $\|\mathcal{R}(X_j)\|_F$ cheaply. 
Our approach gives a whole new  family of algorithmic descriptions of the same approximation sequence $X_j$ to the Riccati solution as each choice of a block basis of $\mathfrak{K}_j^\Box$ leads to a different algorithm.
In case of real $A, B, C$, all iterates are real in case complex shifts are used as complex-conjugate pairs (even so some computations involving complex arithmetic can not be avoided). A new feature of the algorithm, useful for an efficient implementation, is that it allows shifts to be added to the solution not just one at a time, but several at a time. The linear systems of equations of the form $(A^H+\sigma_j)V = C^H$ to be solved for this purpose can be solved simultaneously.
We will rediscover the RADI algorithm \cite{BenBKS18} as a member of the proposed family of algorithms giving not only a new interpretation of the known algorithm but also the new option of adding shifts in parallel. In addition, we specify a member of the proposed family of algorithms which allows for a faster computation of the approximations $X_j$ than any of the three equivalent algorithms from \cite{BenBKS18} in case $m$ is significantly larger than $p$ (and the Sherman-Morrision-Woodbury formula as well as a direct solver is used to solve the linear systems that occur). The results discussed here can be found in slightly different form in the PhD thesis \cite{Ber21}.

We start out with some preliminaries and auxiliary results for the case $p=1$ in Section \ref{sec:p1}. These will be generalized to the case $p>1$ in Section \ref{sec:p>1}. Our main result  is presented in Section \ref{sec:main}. An algorithmic approach suggested by this result is stated. Moreover, we comment on how to adapt the approach so that it can be applied to generalized Riccati equations as well to nonsymmetric Riccati equations. Finally we briefly note that the approximate solution $X_j$ can be interpreted as the solution $X$ of \eqref{ric0} projected onto $\mathfrak{K}_j^\Box(A^H,C^H,\frak{q}_j).$
Section \ref{sec:newalgs} explains how the general algorithmic approach proposed in Section \ref{sec:main} can be implemented in an incremental fashion. This general approach is concretized in Sections \ref{subsec1} and \ref{subsec2} by choosing two specific bases of the underlying block Krylov space. In particular, it is discussed how to add several shifts at a time which allows for a parallel implementation of the time consuming parts of the algorithm. Furthermore, it is discussed how to modify the algorithms in case of real system matrices $A, B, C$ and complex shifts such that the iterates $X_j$ remain real valued. 
Finally, numerical examples comparing the two specific choices of the underlying block Krylov basis are presented in Section \ref{sec_numex}.


\section{Preliminaries and Auxilliary Results for the Case $p=1$}\label{sec:p1}

To keep the notation and the derivations as simple and easy to follow as possible, we will start by considering the case $p=1.$ That is, we consider rational Krylov spaces with a single starting vector $v \in \mathbb{C}^n.$ Rational Krylov spaces were
initially proposed by Ruhe in the 1980s for the purpose of solving large sparse eigenvalue
problems \cite{Ruh84,Ruh94,Ruh98}. In our presentation, we will essentially refer only to the more recent work by G\"uttel and co-authors \cite{BerG15}.
We first recall some definitions and results on rational Krylov subspaces and rational Arnoldi decompositions from \cite[Section 2]{BerG15}. Then we prove some auxiliary results needed in the discussion of our main result presented in Section \ref{sec:main}.

Given a matrix $F\in \mathbb{C}^{n \times n},$ a starting vector $v \in \mathbb{C}^n$, an integer $j <n$ with $\dim\kappa_{j+1}(F,v) = j+1$  and a nonzero polynomial $q_j \in \Pi_j$ with roots disjoint from the spectrum $\Lambda(F),$ a rational Krylov subspace is defined as (see, e.g., \cite[(2.1)]{BerG15}) 
\begin{equation}\label{RK}
\kappa_{j+1}(F,v,q_j) = q_j(F)^{-1}\kappa_{j+1}(F,v) = \{q_j(F)^{-1} p(F)v \mid p \in \Pi_{j}\}
\end{equation}
where $\kappa_{j+1}(F,v) = \{v, Fv, \ldots, F^{j}v\}$ denotes the standard Krylov subspace.  Here $\Pi_j$ denotes the set of polynomials of degree at most $j.$
The roots of $q_j$ are called poles of the rational Krylov space and are denoted by $s_1, \ldots, s_j.$ If the degree $d $ of $q_j$ is less than $j,$  then $j-d$ of the poles are set to $\infty.$ 
The spaces $\kappa_{j+1}(F,v,q_j)$ and $\kappa_{j+1}(F,v)$ are of the same dimension for all $j.$  The poles  of a
rational Krylov space are uniquely determined by the starting vector and vice versa (see \cite[Lemma 2.1]{BerG15}).

There is a one-to-one correspondence between rational Krylov spaces and so-called rational Arnoldi decompositions \cite{Ruh98,BerG15}.
A relation of the form
\begin{equation}\label{rad}
FV_{j+1}\underline{K}_j = V_{j+1} \underline{H}_j
\end{equation}
is called a rational Arnoldi decomposition (RAD) if $V_{j+1} \in \mathbb{C}^{n \times j+1}$ is of full column rank,
$ \underline{H}_j, \underline{K}_j$ are upper Hessenberg matrices of size $(j+1)\times j$ with $|h_{i+1,i}| +|k_{i+1,i}|\neq 0$ for all $i = 1, \ldots, j,$ and the quotients $\mu_i = h_{i+1,i}/k_{i+1,i}$, called poles of the decomposition, are outside the spectrum $\Lambda(F),$ for $i = 1, \ldots, j$ (see \cite[Definition 2.3]{BerG15}).
The columns of $V_{j+1}$ are called the basis of the decomposition and they span the space of the decomposition. 
As noted in \cite{BerG15}, both $\underline{H}_j$ and $\underline{K}_j$ in the RAD \eqref{rad} are of full rank. 

\begin{theorem} \cite[Theorem 2.5]{BerG15} \label{theo1}
Let $\mathcal{V}_{j+1}$ be a vector space of dimension $j+1$. Then $\mathcal{V}_{j+1}$ is a
rational Krylov space with starting vector $v \in \mathcal{V}_{j+1}$ and poles $\mu_1, \ldots, \mu_j$  if and only
if there exists an RAD \eqref{rad} with $\operatorname{range}(V_{j+1}) = \mathcal{V}_{j+1}$, $v_1 = v,$ and poles $s_0=\mu_1,$ $\ldots,$ $s_j=\mu_j$.
\end{theorem}

In our work we will make use of the following two special rational Krylov subspaces. Let $j <n$ be an integer with $\dim \kappa_{j}(F,v) = j.$ Let $\mathfrak{q}_j$ be the polynomial of degree $j$ with roots $s_1, \ldots, s_j \in \mathbb{C}$ and let $\mathcal{S}_j$ denote the set $\{s_1, \ldots, s_j\}$ of roots of $\mathfrak{q}_j.$ Then 
\begin{equation}\label{eq:myK}
\mathfrak{K}_j = \mathfrak{K}_j(F,v,\mathfrak{q}_j) = \mathfrak{q}_j(F)^{-1}\kappa_{j}(F,v) = \{\mathfrak{q}_j(F)^{-1} p(F)v \mid p \in \Pi_{j-1}\}
\end{equation}
defines a rational Krylov subspace.  In order to emphasize the poles of the rational Krylov subspace, we will use the notation $\mathfrak{K}_j(F,v,\mathcal{S}_j)$ instead of $\mathfrak{K}_j(F,v,\mathfrak{q_j}).$ 
The difference to the definition in \eqref{RK} is the choice of the Krylov space, which here is of dimension $j$ while in \eqref{RK} it is of dimension $j+1.$ 
Thus, $\mathfrak{q}_j^{-1}p$ is a proper rational function
and so
\begin{equation}\label{vnotin}
v \notin \mathfrak{K}_j(F,v,\mathfrak{q}_j) .
\end{equation}
In contrast, the rational function $q_j^{-1}p$ in \eqref{RK} may be improper (as $q_j$ may have the same degree as $p$), but can be written as the sum of a polynomial and a proper rational function.

Moreover, we will consider 
\[
\mathfrak{K}_j^+ = \mathfrak{K}_j^+(F,v,\mathcal{S}_j) =\mathfrak{q}_j(F)^{-1}\kappa_{j+1}(F,v) 
\]
 which can be understood as $ \mathfrak{K}_{j+1}(F,v,\mathcal{S}_j \cup \infty).$ This is a rational Krylov space as in \eqref{RK}.
The space $\mathfrak{K}_j^+$ is of dimension $j+1$; it contains the $j$ dimensional spaces $\mathfrak{K}_j$ as well as $F\mathfrak{K}_j$,
\[
F\mathfrak{K}_j=
\mathfrak{q}_j(F)^{-1}F \kappa_j(F,v) \subset \mathfrak{q}_j(F)^{-1}\kappa_{j+1}(F,v) 
=\mathfrak{K}_j^+.
\]
Hence, it is possible to encode the effect of multiplication of $\mathfrak{K}_j$ with $F$ in a decomposition of the form
\[
FV_{j+1}\underline{K}_j = V_{j+1}\underline{H}_j
\]
with a matrix $V_{j+1} \in \mathbb{C}^{n \times j+1}$ whose columns span  $\mathfrak{K}_j^+$, matrices $\underline{K}_j, \underline{H}_j\in \mathbb{C}^{j+1 \times j}$ and 
\begin{align}
\begin{split}\label{eqspans}
\operatorname{span}(V_{j+1}\underline{K}_j)& = \mathfrak{K}_j,\\
\operatorname{span}(V_{j+1}\underline{H}_j)& = F\mathfrak{K}_j.
\end{split}
\end{align}

Now we are ready to prove our first auxiliary result which guarantees the existence of a RAD in a very special form.
\begin{lemma}\label{lemma1}
Let $j <n$ be an integer with $\dim \kappa_{j}(F,v) = j.$ Let $\mathfrak{q}_j$ be a polynomial of degree $j$ with roots $s_1, \ldots, s_j \in \mathbb{C}$ disjoint from the spectrum $\Lambda(F).$  Let $\mathcal{S}_j$ denote the set $\{s_1, \ldots, s_j\}$ of roots of $\mathfrak{q}_j.$ 
Then there exists a RAD
\[
F V_{j+1}\underline{K}_j = V_{j+1} \underline{H}_j
\]
with $\underline{K}_j = \left[\begin{smallmatrix} 0\\I_j \end{smallmatrix}\right] \in \mathbb{R}^{j+1 \times j},$ $\underline{H}_j \in \mathbb{C}^{j+1 \times j}$ is an upper Hessenberg matrix and $V_{j+1}=[v~~Z_j]$ where $\operatorname{span}(Z_j) = \mathfrak{K}_j(F,v,\mathcal{S}_j).$
\end{lemma}
\begin{proof}
As $\dim \kappa_{j}(F,v) =j,$ we have $\dim \mathfrak{K}_j^+ = j+1.$ Thus Theorem \ref{theo1} guarantees the existence of a RAD
\[
F\tilde V_{j+1}\underline{\tilde K}_j = \tilde V_{j+1} \underline{\tilde H}_j
\]
with starting vector $v = \alpha \tilde V_{j+1}e_1, \alpha \in \mathbb{C}.$ 

Moreover, we have from \eqref{vnotin} and \eqref{eqspans} that  $v \not \in \mathfrak{K}_j = \operatorname{span}(\tilde V_{j+1}\underline{K}_j).$ Thus, $e_1 \not \in \operatorname{span}(\underline{K}_j)$. As $\underline{\tilde K}_j$ is a full rank upper Hessenberg matrix, we see that the matrix
$R = [\alpha e_1~~\underline{\tilde K}_j] \in \mathbb{C}^{j+1 \times j+1}$  is nonsingular upper triangular matrix. Thus, 
$\underline{H}_j := R^{-1}\underline{\tilde H}_j$ is an upper Hessenberg matrix,
\[
\underline{K}_j := R^{-1}\underline{\tilde K}_j = \left[\begin{smallmatrix} 0\\I_j \end{smallmatrix}\right] 
\]
and $V_{j+1}:= \tilde V_{j+1}R$ with $V_{j+1} e_1 = v.$ Hence,  $V_{j+1} = [v~~ Z_j]$ and from \eqref{eqspans} we have
\[
\operatorname{span}(V_{j+1}\underline{K}_j) = \operatorname{span}([v~~Z_j] \left[\begin{smallmatrix} 0\\I_j \end{smallmatrix}\right] ) = \operatorname{span}(Z_j) = \mathfrak{K}_j.
\]

\end{proof}

Thus, under the assumptions of Lemma \ref{lemma1} we have with $\underline{H}_j = \left[\begin{smallmatrix}h_j\\ \underline{H}_{-j}\end{smallmatrix}\right],$ where $\underline{H}_{-j} \in \mathbb{C}^{j \times j}$ is an upper triangular matrix, a RAD of the form
\begin{equation}\label{myrad}
FV_{j+1}\left[\begin{smallmatrix} 0\\I_j \end{smallmatrix}\right] = V_{j+1} \left[\begin{smallmatrix}h_j\\ \underline{H}_{-j}\end{smallmatrix}\right].
\end{equation}
Its poles $s_i= \mu_i = h_{i+1,1}/k_{i+1,i} =  h_{i+1,1}, i = 1, \ldots, j$ are just the eigenvalues (that is, the diagonal elements) of $\underline{H}_{-j}.$

With the additional assumption $\mathcal{S}_j \cap -\overline{\mathcal{S}_j} = \emptyset$ we have that $\underline{H}_{-j}$ and $-\underline{H}_{-j}^H$ have no eigenvalues in common. Hence, for any $W \in \mathbb{C}^{j \times m}$ the Lyapunov equation
\begin{equation}\label{lyap}
\tilde Y_j \underline{H}_{-j} + \underline{H}_{-j}^H \tilde Y_j - WW^H= 0
\end{equation}
has a unique Hermitian solution $\tilde Y_j$ (see, e.g., \cite[Theorem 4.4.6]{HorJ94} or \cite[Theorem 5.2.2]{LanR95}). The assumption $\mathcal{S}_j \cap -\overline{\mathcal{S}_j} = \emptyset$ holds in particular if all $s_i \in \mathcal{S}_j$ lie in the open right half plane (that is, $\operatorname{Re}(s_i) > 0, i = 1, \ldots, j$). In that case the following lemma holds.

\begin{lemma}\label{lemma2}
In addition to  the assumptions of Lemma \ref{lemma1} suppose that all $s_i \in \mathcal{S}_j$ lie in the open right half plane.
Moreover, let $WW^H= Z_j^HBB^HZ_j +h_j^Hh_j$ for some $B \in \mathbb{C}^{n \times m}.$
Then there exists a unique positive definite solution of the Lyapunov equation \eqref{lyap}.
\end{lemma}
\begin{proof}
As $WW^H$ is positive semi-definite, we have from \cite[Theorem 5.3.1 (a)]  {LanR95} that $\tilde Y_j$ is positive semi-definite.

In case  $(\underline{H}_{-j}^H, h_j^H)$ is controllable (that is, $\operatorname{rank}([h_j^H~~ \underline{H}_{-j}^H h_j^H~~ \cdots~~(\underline{H}_{-j}^H)^{j-1} h_j^H] )=j$), the pair  $(\underline{H}_{-j}^H, [h_j^H~~Z_j^HB])$ is controllable as well. Thus, with  \cite[Theorem 5.3.1 (b)] {LanR95} we obtain that $\tilde Y_j$ is positive definite.

Now assume to the contrary that $(\underline{H}_{-j}^H, h_j^H)$ is not controllable. This implies that $(h_j, \underline{H}_{-j})$ is not observable,  that is, $\operatorname{rank}([h_j^H~~ \underline{H}_{-j}^H h_j^H~~ \cdots~~(\underline{H}_{-j}^H)^{j-1} h_j^H]^H )<j$ see, e.g., \cite[Theorem 4.2.2]{LanR95}). Then, for any right eigenvector $z \in \mathbb{C}^j\backslash \{0\}$ of $\underline{H}_{-j}$ we have $h_j^Hz =0$ (see, e.g, \cite[Theorem 4.3.3]{LanR95} or \cite[Lemma 3.3.7]{Son98}). As the eigenvalues of $\underline{H}_{-j}$ are just the $s_i$, it holds $\underline{H}_{-j} z= s_i z$ for some $i.$ Moreover,
\[
FV_{j+1}\underline{H}_{-j} z = FV_{j+1}\left[\begin{smallmatrix}h_j\\\underline{H}_{-j} \end{smallmatrix}\right]z =s_i  FV_{j+1}\left[\begin{smallmatrix}0\\ I_{j} \end{smallmatrix}\right]z = s_i V_{j+1}\underline{H}_{-j}z
\]
where the last equality is due to Lemma \ref{lemma1}. Hence, $V_{j+1}\underline{H}_{-j} z$ is an eigenvector of $F$ with eigenvalue $s_i.$ This is a contradiction to the definition of a RAD because the poles $\mu$ must be distinct from the eigenvalues of $F$.
\end{proof}
Hence, once the basis of the Krylov subspace $\mathfrak{K}_j(F,v,\mathcal{S}_j)$ is fixed (via the columns of $Z_j$, see Lemma \ref{lemma1}),  the Lyapunov equation \eqref{lyap} has a unique solution.

\begin{corollary}
Let the assumptions of Lemma \ref{lemma2} hold. Then there exists a unique positive-definite solution of the Riccati equation
\begin{equation}\label{ric1}
\underline{H}_{-j} Y_j + Y_j \underline{H}_{-j}^H - Y_j(Z_j^HBB^HZ_j+h_j^Hh_j)Y_j = 0.
\end{equation}
\end{corollary}
\begin{proof}
Due to Lemma \ref{lemma2} there exists  a unique positive definite solution of the Lyapunov equation $\tilde Y_j \underline{H}_{-j} + \underline{H}_{-j}^H \tilde Y_j - WW^H= 0$. With $Y_j = \tilde Y_j^{-1}$ this a equivalent to \eqref{ric1}.
\end{proof}
\section{Preliminaries and Auxiliary Results for $p > 1$}\label{sec:p>1}

In this section, the results presented in the previous section are generalized to the case $p>1.$ This implies that instead of a single starting vector $v \in \mathbb{C}$, we now have to deal with a block starting vector $\textbf{v} \in \mathbb{C}^{n \times p}.$ First some definitions and results on block rational  Krylov subspaces and block rational Arnoldi decompositions from \cite[Section 1 and 2]{ElsG20} are recalled. Then we explain how to generalize the auxiliary results obtained in the previous section to  case of a block starting vector.

Given a matrix $F\in \mathbb{C}^{n \times n}$ and a starting block vector $\textbf{v} \in \mathbb{C}^{n\times p}$ of maximal rank, the associated block Krylov subspace of order $j+1$ is defined as
\begin{align*}
\kappa_{j+1}^\Box(F,\textbf{v}) &= \operatorname{blockspan}\{\textbf{v}, F\textbf{v}, \ldots, F^{j}\textbf{v}\} = \{ \sum_{k=0}^{j}F^k\textbf{v} \frak{F}_k \mid \frak{F}_k \in \mathbb{C}^{p \times p}\}.
\end{align*}
From here on it is assumed that the $(j+1)p$  columns of $\kappa_{j+1}^\Box(F,\textbf{v})$ are linearly independent. Then, the block Krylov subspace  $\kappa_{j+1}^\Box(F,\textbf{v})$ has dimension $(j+1)p^2,$ and
 every block vector $\sum_{k=0}^{j}F^k\textbf{v} \frak{F}_k$ $\in  \kappa_{j+1}^\Box(F,\textbf{v})$ corresponds to exactly one matrix polynomial $P(z)=\sum_{k=0}^{j}z^k\frak{F}_k,$
\[
P(F)\circ \textbf{v} = \textbf{v} \frak{F}_0 + F\textbf{v}\frak{F}_1 + F^2\textbf{v}\frak{F}_2 + \cdots + F^{j} \textbf{v} \frak{F}_j,
\]
see, e.g., \cite{ElsG20}.
Given a nonzero polynomial $q_j \in \Pi_j$ with roots disjoint from the spectrum $\Lambda(F),$ a block rational Krylov subspace is defined as 
\begin{equation}\label{RBK}
\kappa_{j+1}^\Box(F,\textbf{v},q_j) = q_j(F)^{-1}\kappa_{j+1}^\Box(F,\textbf{v}) .
\end{equation}
The roots of $q_j$ are called poles of the block rational Krylov space and are denoted by $s_1, \ldots, s_j.$ 
The spaces $\kappa_{j+1}^\Box(F,\textbf{v},q_j)$ and $\kappa_{j+1}^\Box(F,\textbf{v})$ are of the same dimension for all $j.$  

There is a one-to-one correspondence between block rational Krylov spaces and so-called block rational Arnoldi decompositions (BRAD) (that is, an analogue of Theorem \ref{theo1} holds here). 
A relation of the form
\begin{equation}\label{brad}
FV_{j+1}\underline{K}_j = V_{j+1} \underline{H}_j
\end{equation}
is called a block rational Arnoldi decomposition (BRAD) if $V_{j+1} \in \mathbb{C}^{n \times (j+1)p}$ is of full column rank,
$ \underline{H}_j, \underline{K}_j$ are block upper Hessenberg matrices of size $(j+1)p\times jp$ where at least one of the matrices $H_{i+1,i}$ and $K_{i+1,i}$ is nonsingular, $\beta_iK_{i+1,i} =\gamma_iH_{i+1,i}$ with scalars $\beta_i,\gamma_i \in \mathbb{C}$ such that $|\beta_i|+|\gamma_i|\neq 0$ for all $i = 1, \ldots, j,$ and the quotients $\mu_i = \beta_i/\gamma_i$, called poles of the BRAD, are outside the spectrum $\Lambda(F),$ for $i = 1, \ldots, j$ (see \cite[Definition 2.2]{ElsG20}).
As noted in \cite[Lemma 3.2 (ii)]{ElsG20}, in case one of the subdiagonal blocks $H_{i+1,i}$ and $K_{i+1,i}$ is singular, it is the zero matrix. 
The block columns of $V_{j+1}=[v_1 ~\ldots~v_{j+1}]$ blockspan the space of the decomposition, that is, the linear space of block vectors $v = \sum_{k=1}^{j+1}v_k\frak{F}_k$ with arbitrary coefficient matrices $\frak{F}_k \in \mathbb{C}^{p \times p}.$ An algorithm which constructs a BRAD can be found in \cite[Algorithm 2.1]{ElsG20}, see also \cite{BerEG14}.

In analogy to the rational Krylov space $\mathfrak{K}_j(F,v,\mathfrak{q}_j) $ \eqref{eq:myK} we will make use of the following special block rational Krylov subspace
\[
\mathfrak{K}_j^\Box = \mathfrak{K}_j^\Box(F,\textbf{v},\mathcal{S}_j)= \mathfrak{K}_j^\Box(F,\textbf{v},\mathfrak{q}_j) = \mathfrak{q}_j(F)^{-1}\kappa_{j}^\Box(F,\textbf{v}) .
\]
Here $j$ is an integer with $\dim \kappa_{j}^\Box(F,\textbf{v}) = jp^2$ and $\mathfrak{q}_j$ is the polynomial of degree $j$ with roots $s_1, \ldots, s_j \in \mathbb{C}.$ As before, $\mathcal{S}_j$ denotes the set $\{s_1, \ldots, s_j\}$ of roots of $\mathfrak{q}_j.$ 
The difference to the definition in \eqref{RBK} is the choice of the Krylov space, which here is of order $j$ and dimension $jp^2,$ while in \eqref{RBK} it is of order $j+1$ and dimension $(j+1)p^2.$ 

Now we can generalize the auxiliary results from the previous section to the case $p>1.$ Essentially one needs to replace vectors with block vectors and scalars by scalar $p \times p$ matrices (i.e. multiples of $I_p$). Let a set $\mathcal{S}_j$ of $j$ roots of a polynomial $\mathfrak{q}_j$ be given with $\mathcal{S}_j \cap \Lambda(F) = \emptyset.$ The generalization of Lemma \ref{lemma1} yields a BRAD
\begin{equation}\label{brad1}
FV_{j+1}\underline{K}_j = V_{j+1}\underline{H}_j
\end{equation}
with $V_{j+1} = [ \textbf{v}~~ Z_j] \in \mathbb{C}^{n \times (j+1)p}$ where $\operatorname{blockspan}(Z_j) = \mathfrak{K}^\Box_j(F,\textbf{v},\mathcal{S}_j)$ and block upper Hessenberg matrices
\begin{align*}
\underline{K}_j &= \begin{bmatrix} 0 \\ I_{jp}\end{bmatrix} \in \mathbb{R}^{(j+1)p\times jp}\\
\underline{H}_j &= \begin{bmatrix} h_j\\ \underline{H}_{-j}\end{bmatrix} =\begin{bmatrix}
 H_{11} & H_{12} &\cdots & H_{1,p-1}& H_{1p}\\
\mu_1I_p & H_{22} &\cdots & H_{2,p-1}& H_{2p}\\
 & \ddots & \ddots & \vdots & \vdots\\
&&\mu_{j-2}I_p&H_{j-1,p-1}& H_{j-1,p}\\
&&&\mu_{j-1}I_p& H_{jp}\\
&&&& \mu_{j}I_p
\end{bmatrix}\in \mathbb{C}^{(j+1)p\times jp} 
\end{align*}
with $h_j \in \mathbb{C}^{p \times jp},$ $\underline{H}_{-j} \in \mathbb{C}^{jp \times jp}$ and $H_{ij}\in \mathbb{C}^{p \times p}.$

Moreover, as in Lemma \ref{lemma2} we have that the Lyapunov equation
\begin{equation}\label{lyap_block}
\tilde Y_j \underline{H}_{-j} + \underline{H}_{-j}^H \tilde Y_j - Z_j^HBB^HZ_j-h_j^Hh_j = 0
\end{equation}
 for some $B \in \mathbb{C}^{n \times m}$ has a unique positive definite solution.

\section{Main Result}\label{sec:main}
Now we turn our attention to solving the continuous-time algebraic Riccati equation
\begin{equation}\label{ric}
0 = \mathcal{R}(X) := A^HX+XA+C^HC-XBB^HX 
\end{equation}
with $A \in \mathbb{C}^{n \times n}, B\in \mathbb{C}^{n \times m}$ and $C\in \mathbb{C}^{p \times m}.$ 

For our discussion in this section, we choose a fixed set of $j$ roots $\mathcal{S}_j \subset \mathbb{C}_+$  with $\mathcal{S}_j \cap \Lambda(A) = \emptyset$  and the corresponding block rational Krylov subspace $\mathfrak{K}_j^\Box=\mathfrak{K}_j^\Box(A^H,C^H,\mathcal{S}_j).$ Moreover, we assume that the $jp$ columns in $\kappa_j^\Box(A^H,C^H)$ are linearly independent. That is, the assumptions of the previous section hold. We are interested in an approximate solution $X_j$ of \eqref{ric} which satisfies
\begin{itemize}
\item $X_j$ is of rank $jp$ and of the form $X_j = Z_jY_jZ_j^H,$  where
\item $Z_j \in \mathbb{C}^{n \times jp}$  such that $\operatorname{blockspan}(Z_j) =\mathfrak{K}_j^\Box,$  and
\item $\operatorname{rank}(\mathcal{R}(X_j)) = p.$ 
\end{itemize}
We will see that such an $X_j$ is unique. Fixing the initial guess as $X_0=0,$  these iterates will be equal to those in \eqref{eq:equiv1} in case the same shifts are used.

Let \eqref{brad1} hold with $F = A^H$ and $\textbf{v} = C^H.$ Let $X_j = Z_jY_jZ_j^H$ for some Hermitian matrix $Y_j \in \mathbb{C}^{jp \times jp}.$ Then, as $Z_j =V_{j+1}\left[\begin{smallmatrix}0\\I_{jp}\end{smallmatrix}\right] =  V_{j+1}\underline{K}_j$ we can write $X_j = V_{j+1}\underline{K}_jY_j\underline{K}_j^HV_{j+1}^H$ and
\begin{align*}
\mathcal{R}(X_j) &= A^HX_j+X_jA+C^HC-X_jBB^HX_j \\
&= A^H V_{j+1}\underline{K}_jY_j \underline{K}_j^HV_{j+1}^H+ V_{j+1}\underline{K}_jY_j \underline{K}_j^HV_{j+1}^HA+C^HC
\\&\qquad
- V_{j+1}\underline{K}_jY_j \underline{K}_j^HV_{j+1}^HBB^H V_{j+1}\underline{K}_jY_j \underline{K}_j^HV_{j+1}^H.
\end{align*}
With \eqref{brad1} and $S_j = \underline{K}_j^HV_{j+1}^HBB^H V_{j+1}\underline{K}_j = Z_j^HBB^HZ_j$ it follows
\begin{align}
\mathcal{R}(X_j) 
&= V_{j+1}\underline{H}_jY_j \underline{K}_j^HV_{j+1}^H+ V_{j+1}\underline{K}_jY_j \underline{H}_j^HV_{j+1}^H+C^HC- V_{j+1}\underline{K}_jY_j S_jY_j \underline{K}_j^HV_{j+1}^H \nonumber\\
&= V_{j+1}\left( \underline{H}_jY_j \underline{K}_j^H+ \underline{K}_jY_j \underline{H}_j^H+\left[\begin{smallmatrix} I_p\\0\end{smallmatrix}\right] \left[\begin{smallmatrix} I_p\\0\end{smallmatrix}\right]^H - \underline{K}_jY_j S_jY_j \underline{K}_j^H\right)V_{j+1}^H\nonumber\\
&=: V_{j+1}M_jV_{j+1}^H\label{eq1}
\end{align}
as $C^H = V_{j+1}\left[\begin{smallmatrix} I_p\\0\end{smallmatrix}\right].$ 
As $V_{j+1}$ is of full rank, the rank of $\mathcal{R}(X_j)$ is the same as the rank of $M_j.$ We rewrite $M_j$ as
\begin{align}
M_j &= \underline{H}_jY_j \underline{K}_j^H+ \underline{K}_jY_j \underline{H}_j^H+\left[\begin{smallmatrix} I_p\\0\end{smallmatrix}\right] \left[\begin{smallmatrix} I_p\\0\end{smallmatrix}\right]^H - \underline{K}_jY_j S_jY_j \underline{K}_j^H \nonumber\\
&= \begin{bmatrix}0 & h_j Y_j\\ 0& \underline{H}_{j-1}Y_j\end{bmatrix} + \begin{bmatrix} 0 & 0\\ Y_jh_j^H & Y_j \underline{H}_{-j}^H\end{bmatrix}+\begin{bmatrix} I_p & 0\\ 0 & 0\end{bmatrix}
- \begin{bmatrix} 0 & 0\\ 0 & Y_jS_jY_j\end{bmatrix} \nonumber\\
&= \begin{bmatrix} I_p & &h_jY_j\\ Y_jh_j^H  & ~&\underline{H}_{j-1}Y_j + Y_j \underline{H}_{-j}^H -Y_jS_jY_j \end{bmatrix}.\label{eq_Mj}
\end{align}
In case $Y_j$ solves the block version of \eqref{ric1} it follows
\begin{align*}
M_j &=  \begin{bmatrix} I_p & &h_jY_j\\ Y_jh_j^H  & ~&Y_jh_j^Hh_jY_j \end{bmatrix}= \begin{bmatrix} I_p \\ Y_jh_j^H  \end{bmatrix} \begin{bmatrix} I_p & h_jY_j\end{bmatrix}.
\end{align*}
Then, by construction, $M_j$ is a matrix of rank $p.$ Thus, the residual \eqref{eq1} is of rank $p.$ We have
\begin{equation}\label{eq:Rj}
\mathcal{R}(X_j)  = R_jR_j^H
\end{equation}
with $R_j = V_{j+1}\left[\begin{smallmatrix}I_p\\Y_jh_j^H\end{smallmatrix} \right]= C^H + Z_jY_jh_j^H.$ This implies $R_0=C^H$ and $X_0=0.$

This finding is summarized in the next theorem.
\begin{theorem}\label{theo:main}
Let $j <n$ be an integer such that the $jp$ columns of  $\kappa_{j}^\Box(A^H,C^H) $ are linearly independent. Let $\mathfrak{q}_j$ be a polynomial of degree $j$ with roots $s_1, \ldots, s_j \in \mathbb{C}_+$ and $\mathcal{S}_j \cap \Lambda(A) = \emptyset$ where $\mathcal{S}_j$ denotes the set $\{s_1, \ldots, s_j\}$ such that \eqref{brad1} holds. Denote the unique positive-definite solution of 
\[
\underline{H}_{-j} Y_j + Y_j \underline{H}_{-j}^H - Y_j(Z_j^HBB^HZ_j+h_j^Hh_j)Y_j = 0
\]
  by $Y_j$. Then $X_j = Z_jY_jZ_j^H$  is the unique matrix of rank $jp$ such that the residual $\mathcal{R}(X_j)$ is of rank $p.$
\end{theorem}
\begin{proof}
It remains to prove that $X_j$ is unique, that is, there is no other matrix $Y_j$ such that the residual $\mathcal{R}(Z_jY_jZ_j^H)$ is of rank $p.$ In order to see this, let us assume that $M_j$ in \eqref{eq_Mj} is of rank $p.$ Then the first $p$ rows and columns of $M_j$ in \eqref{eq_Mj} determine the rank-$p$-factorization of $M_j$ uniquely,
\[
M_j = \begin{bmatrix} I_p \\ Y_jh_j^H  \end{bmatrix} \begin{bmatrix} I_p & h_jY_j\end{bmatrix}.
\]
Thus $Y_j$ must be the solution of \eqref{ric1}.
\end{proof}

\begin{remark}
There are other $X_j = Z_jY_jZ_j^H$ which yield a rank-$p$-residual $\mathcal{R}(X_j)$, but these require that $X_j$ is no longer of rank $jp.$ An easy example is the choice $Y_j=0.$ This gives the rank-$p$-residual $\mathcal{R}(X_j) = C^HC.$ 
\end{remark}

Recall from the Introduction that the four Riccati ADI methods 
for solving the Riccati equation \eqref{ric} produce (at least theoretically) the same approximate solution $X_j$ of rank $jp$ with $\operatorname{rank}(\mathcal{R}(X_j)) = p.$   Moreover, if the block columns of $Z_j \in \mathbb{C}^{n \times jp}$ represent a basis for the block rational Krylov subspace $\mathcal{L}_j^\Box = \mathcal{L}_j^\Box(A^H,C^H,\Sigma_j)$ \eqref{eq_L}, then $X_j= Z_jY_jZ_j^H$ holds for some matrix $Y_j \in \mathbb{C}^{jp \times jp}$  \cite[Proposition 2]{BenBKS18}.
As for $\Sigma_j = - \mathcal{S}_j,$ we have  $\mathcal{L}_j^\Box  =\mathfrak{K}_j^\Box$ and as $X_j$ is unique, $X_j = X_j^\text{qadi} = X_j^\text{cay} = X_j^\text{radi}$ needs to hold for $X_j$ as in Theorem \ref{theo:main}. 
That is, Theorem \ref{theo:main} implies equivalence of all methods which
yield an approximate solution $X_j$ with a rank-$p$ residual and whose block columns blockspan a Krylov subspace $\mathfrak{K}_j^\Box.$ In particular, this includes the four Riccati ADI methods. Due to the structure of the approximate solution $X_j$ all choices of a basis $Z_j$ of the block Krylov subspace are equivalent, as a transition matrix for a change of basis can be
incorporated into $Y_j.$

Theorem \ref{theo:main} suggests the algorithmic approach summarized in Algorithm \ref{alg1} for solving the algebraic Riccati equation \eqref{ric}. Our point of view is analytical rather than numerical, so no form of orthogonality of $V_{j+1}$ is enforced.
We will see in the next section how the approximate solution $X_j$ can be computed recursively avoiding the explicit upfront construction of the BRAD in Step 1 and the explicit solution of the Lyapunov equation in Step 2. 

\begin{algorithm}
\caption{Algorithmic Approach suggested by Theorem \ref{theo:main}}\label{alg1}
\begin{algorithmic}[1]
\Require $A \in \mathbb{C}^{n \times n}, B\in \mathbb{C}^{n \times m}, C\in \mathbb{C}^{p \times n},$ and a set of shifts $\mathcal{S}_j = \{s_1, \ldots, s_j\}, s_i \in \mathbb{C}_+$ with $\mathcal{S}_j \cap \Lambda(A) = \emptyset.$
\Ensure approximate solution $X_j$ of \eqref{ric}, residual factor $R_j$ such that $\mathcal{R}(X_j)=R_jR_j^H.$
\State Construct BRAD $A^HV_{j+1}\underline{K}_j = V_{j+1}\underline{H}_j$ with $V_{j+1} = [ C^H~~ Z_j],$ $\underline{K}_j = \left[\begin{smallmatrix}0\\I_{jp}\end{smallmatrix}\right],$ and $\underline{H}_j = \left[\begin{smallmatrix} h_j\\ \underline{H}_{-j}\end{smallmatrix}\right].$
\State Solve $\tilde Y_j \underline{H}_{-j} + \underline{H}_{-j}^H \tilde Y_j - Z_j^HBB^HZ_j-h_jh_j^H= 0$ for $\tilde Y_j=\tilde Y_j^H.$
\State Set $X_j = Z_j\tilde Y_j^{-1}Z_j^H$.
\State Set $\tilde h = \tilde Y_j^{-1}h_j^H.$
\State $R_j = V_{j+1}\left[ \begin{smallmatrix} I_p\\\tilde h\end{smallmatrix}\right].$  \Comment{$\mathcal{R}(X_j) = V_{j+1}\left[ \begin{smallmatrix} I_p\\    \tilde Y_j^{-1}h_j^H \end{smallmatrix}\right] \left[ \begin{smallmatrix} I_p~~   \tilde Y_j^{-1}h_j^H \end{smallmatrix}\right] V_{j+1}^H=R_jR_j^H$.}
\end{algorithmic}
\end{algorithm}

\begin{remark}\label{rem_lyap}
Algorithm \ref{alg1} can be applied to continuous Lyapunov equations, 
\begin{equation}\label{eq_lyapsolv}
0 = A^HX + XA + C^HC,
\end{equation}
as these are a special case of the Riccati equation \eqref{ric0} (with $B =0$).  
\end{remark}
\subsection{Generalized Riccati equations}\label{subsec:genR}
Our approach can be adapted for solving the generalized Riccati equation
\begin{equation}\label{ric_gen}
A^HXE+E^HXA+C^HC-E^HXBB^HXE = 0
\end{equation}
with an additional nonsingular matrix $E\in \mathbb{C}^{n \times n}.$ As noted in \cite[Section 4.4]{BenBKS18}, the equivalent Riccati equation
\[
\mathcal{R}_{gen}(X) = E^{-H}A^HX+XAE^{-1}+E^{-H}C^HCE^{-1}-XBB^HB = 0
\]
has the same structure as \eqref{ric0} where the system matrix $A$ and the initial residual factor $C^H$ are replaced by$ AE^{-1}$ and $E^{-H}C^H$, respectively.  In an efficient iteration inverting $E$ is avoided by utilizing the relation $((AE^{-1})^H -\mu I_n)^{-1}E^{-H}R = (A^H -\mu E^H)^{-1}R$
with a residual factor $R$ of \eqref{ric_gen}. This requires some (standard) modifications in the algorithms presented.

Our approach can also be applied to the nonsymmetric Riccati equation
\begin{equation}\label{ric_nonsym}
\mathcal{R}_{nonsym}(X) = A_1^HX+XA_2+C_1^HC_2-XB_2B_1^HX = 0
\end{equation}
with $A_i \in \mathbb{C}^{n_i \times n_i}, B_i \in \mathbb{C}^{n_i \times m},$ and $C_i \in \mathbb{C}^{p \times n_i}$ for $i = 1, 2$ and $X \in \mathbb{C}^{n_1\times n_2}.$ Consider the two decompositions
\[
A_i^H\left[C_i^H~~Z_i\right] \underline{K}_i = \left[C_i^H~~Z_i\right] \underline{H}_i
\]
for $i = 1, 2,$ with $\underline{K}_i = \left[\begin{smallmatrix}0\ I\end{smallmatrix}\right]$ and $\underline{H}_i = \left[\begin{smallmatrix}h_i\\\underline{H}_{-i}\end{smallmatrix}\right]$  where the number of columns of $Z_1$ and $Z_2$ are the same. Then in analogy to \eqref{eq1} we can rewrite the residual \eqref{ric_nonsym} for the solution $X = Z_1YZ_2^H$ as
\[
\left[C^H_1~~Z_1\right] \left( \underline{H}_1Y\underline{K}_2^H+\underline{K}_1Y\underline{H}_2^H +
\begin{bmatrix}I_p\\0\end{bmatrix} \left[ I_p~~0\right] - \underline{K}_1YSY\underline{K}_2^H\right)
\left[C_2^H~~Z_2\right]^H
\]
due to $Z_i=\left[C_i^H~~Z_i\right]\underline{K}_i$ and with $S=Z_2^HB_2B_1^HZ_1.$  $Y^{-1}$ is no longer determined by the Lyapunov equation \eqref{lyap_block}, it is now determined by the Sylvester equation
\[
0=Y^{-1}\underline{H}_{-1}+\underline{H}_{-2}^HY^{-1}+(S+h_2^Hh_1),
\]
All algorithms presented in the following can be adapted to take care of a nonsymmetric Riccati equation \eqref{ric_nonsym} by some (more or less straightforward) modifications.

\subsection{Projection} 
The Riccati ADI approximate solution $X_j$ can be interpreted as the solution of a projection of the large-scale Riccati equation \eqref{ric0} onto the Krylov subspace $\mathfrak{K}_j(A^H,C^H,\mathfrak{q}_j)$ \eqref{eq:myK} described by $Z_j.$

We will consider only the case $p=1$ as the final result will hold only for that choice of $p.$ Let
\[
A^HV_{j+1}\begin{bmatrix}0\\I_j\end{bmatrix}=V_{j+1}\begin{bmatrix}h_j\\\underline{H}_{-j}\end{bmatrix}
\]
hold as in Theorem \ref{theo:main} with $V_{j+1}=\left[ C^H~~Z_j\right].$ Let $W\in \mathbb{C}^{n \times j}$ be of rank $j$ such that $Z_j^HW$ is nonsingular. Then $\Pi = Z_j(W^HZ_j)^{-1}W^H \in \mathbb{C}^{n \times n}$ is a projection onto $\text{im}(\Pi)=\text{im}(Z_j)$ along the kernel of $\Pi,$ where $\text{ker}(\Pi)=\text{ker}(W^H).$ Set $\widetilde W = W(Z_j^HW)^{-1}$. Then $\Pi = Z_j\widetilde W$ and $\widetilde W^HZ_j=I.$

Assume further that $\widetilde W^H R_j =0,$ that is, that the residual factor $R_j$ with $\mathcal{R}(X_j)=R_jR_j^H$ lies in the kernel of the projection $\Pi.$  Then $0=\widetilde W^HR_j = \widetilde W^H(C^H+Z_jY_jh_j^H)$ holds, that is, $\widetilde W^HC^H = -Y_jh_j^H.$ Moreover, $\Pi \mathcal{R}(X_j)\Pi^H=0$ must hold. As $\text{ker}(Z_j)=\{0\}$, this is equivalent to
\begin{align*}
0 &= \widetilde W^H\mathcal{R}(X_j)\widetilde W \\
&= \widetilde W^H V_{j+1} 
\begin{bmatrix} 1 & h_jY_j\\ Y_jh_j^H & \underline{H}_{-j}Y_j+Y_j \underline{H}_{-j}^H-Y_jS_jY_j\end{bmatrix}V_{j+1}^H\widetilde W\\
&= -Y_jh_j^Hh_jY_j+\underline{H}_{-j}Y_j+Y_j\underline{H}_{-j}-Y_jS_jY_j,
\end{align*}
where we used \eqref{eq1} and \eqref{eq_Mj} and the observation that $\widetilde W^HV_{j+1} = [-Y_jh_j^H~~I].$
Hence, the projected equation $\Pi \mathcal{R}(X_j) \Pi^H=0$ is equivalent to the small scale Riccati equation stated in Theorem \ref{theo:main}. A similar observation for the ADI iteration to solve Lyapunov equations has been made in \cite[Section 3.2]{WolP16} and \cite[Remark 5.16]{Wol14}. 

Although $\widetilde W$ is unknown in practice, the projected system matrices $\widetilde W^HA^HZ_j$ and $\widetilde W^H(A^H -X_jBB^H)Z_j$ can be expressed in terms of parts of the RAD,
\begin{align*}
\widetilde W^HA^HZ_j &= -Y_jh_j^Hh_j+\underline{H}_{-j},\\
\widetilde W^H(A^H -X_jBB^H)Z_j &= -Y_j\underline{H}_{-j}^HY_j^{-1}.
\end{align*}
Hence, the eigenvalues of the projected matrix $\widetilde W^H(A^H -X_jBB^H)Z_j$ are the negative conjugate poles of the underlying Krylov subspace which are given by the eigenvalues of $\underline{H}_{-j}.$

The residual factor $R_j = V_{j+1}\left[\begin{smallmatrix} 1\\Y_jh_j^H\end{smallmatrix}\right]$ is a linear combination of the columns of $V_{j+1}$ and so it is an element of $\mathfrak{K}_j^+.$ In other words, $R_j$ can be represented as a rational function in $A^H$ multiplied with $C^H.$ In particular, it holds that
\[
R_j = \mathfrak{q}_j(A^H)^{-1}\mathfrak{p}_j(A^H)C^H
\]
with 
\[
\mathfrak{q}_j(x) = \prod_{i=1}^j(x-s_i)\quad \text{and}\quad
\mathfrak{p}_j(x) = \prod_{i=1}^j(x-\lambda_i^{(j)}),
\]
where the $s_j$ are the poles of the Krylov subspace and the $\lambda_i^{(j)}$ are the eigenvalues  of the projected matrix $\widetilde W^HA^HZ_j = -Y_jh_j^Hh_j+\underline{H}_{-j}.$

\section{New Algorithms}\label{sec:newalgs}
In this section, we will first discuss how Algorithm \ref{alg1} can be implemented in an incremental fashion.
In order to achieve this, it is observed that when increasing the order of the block Krylov space by one, then $Y_{j+1}^{-1} = \left[ \begin{smallmatrix} Y_j^{-1} & Y_{12}\\Y_{12}^H & Y_{22}\end{smallmatrix}\right]$ holds where $Y_{12}$ can be computed by solving a linear system and $Y_{22}$ can be computed directly. Thus, there is no need to solve the Lyapunov equation \eqref{lyap}, resp., the Riccati equation \eqref{ric1} associated with $\mathfrak{K}^\Box_{j+1}.$

Let $\mathfrak{K}_{j+1}^\Box(A^H,C^H,\mathcal{S}_{j+1})$ with an associated BRAD be  given as in \eqref{brad1}
\begin{equation}\label{brad_here}
A^H \left[ C^H~~Z_{j+1}\right] \begin{bmatrix} 0\\I_{(j+1)p} \end{bmatrix} = \left[ C^H~~Z_{j+1}\right] \begin{bmatrix}h_{j+1}\\\underline{H}_{-(j+1)}\end{bmatrix}
\end{equation}
and $X_{j+1} = Z_{j+1}Y_{j+1}Z_{j+1}^H$ with $Y_{j+1} = \tilde Y_{j+1}^{-1}$ where $\tilde Y_{j+1}$ is the solution of 
\begin{equation}\label{lyap_jp1}
\tilde Y_{j+1} \underline{H}_{-(j+1)} + \underline{H}_{-(j+1)}^H \tilde Y_{j+1} - Z_{j+1}^HBB^HZ_{j+1}-h_{j+1}^Hh_{j+1} = 0
\end{equation}
This BRAD includes the BRAD associated with $\mathfrak{K}_{j}^\Box(A^H,C^H,\mathcal{S}_{j})$ in its first $jp$ columns as
\begin{equation}\label{eq:brad}
A^H \left[ C^H~~Z_j~~\hat Z\right] \begin{bmatrix} 0 & 0\\I_{jp}&0\\ 0 & I_p \end{bmatrix} = \left[ C^H~~Z_j~~\hat Z\right] \begin{bmatrix}h_j & U_1\\\underline{H}_{-j}&U_2 \\ 0 & D\end{bmatrix},
\end{equation}
that is, 
\begin{align*}
Z_{j+1} &= \left[Z_j~~\hat Z\right], \qquad
h_{j+1} = \left[ h_j ~~ U_1  \right], \qquad
\underline{H}_{-(j+1)} =  \begin{bmatrix} \underline{H}_{-j}&U_2 \\ 0 & D\end{bmatrix},
\end{align*}
where $D = \mu_{j+1}I_p \in \mathbb{C}^{p \times p}, U_1 \in \mathbb{C}^{p\times p}$ and $U_2\in \mathbb{C}^{jp\times p}.$
Hence, we have $X_j=Z_jY_jZ_j^H$ with  $Y_{j} = \tilde Y_j^{-1}$ as in \eqref{lyap_block},
\begin{equation}\label{lyap_here}
\tilde Y_{j} \underline{H}_{-j} + \underline{H}_{-j}^H \tilde Y_{j} - Z_{j}^HBB^HZ_{j}-h_{j}^Hh_{j} = 0.
\end{equation}
The choice of $\hat Z$ in \eqref{eq:brad} describes the basis of $\mathfrak{K}_{j+1}^\Box(A^H,C^H,\mathcal{S}_{j+1})$ used. Each possible $\hat Z$ yields a different solution $Y_{j+1}.$  Recall, that each possible $\hat Z$ in \eqref{eq:brad} provides the same approximate solution $X_{j+1}$.

We will see next that there is no need to solve \eqref{lyap_jp1} for $\tilde Y_{j+1}$ as $\tilde Y_{j+1}$ can be obtained from $\tilde Y_j.$  
Partition the solution $Y_{j+1}^{-1} = \tilde Y_{j+1}$ of \eqref{lyap_jp1}  as
\[
Y_{j+1}^{-1} = \tilde Y_{j+1} = \begin{bmatrix} Y_{11} & Y_{12}\\Y_{12}^H & Y_{22}\end{bmatrix}
\]
with $Y_{11}\in \mathbb{C}^{jp \times jp}, Y_{12} \in  \mathbb{C}^{jp \times p}$ and $Y_{22} \in  \mathbb{C}^{p \times p}.$
Then \eqref{lyap_jp1} reads with $D = \mu_{i+1}I_p$
\begin{align*}
0 &= \begin{bmatrix} Y_{11} & Y_{12}\\Y_{12}^H & Y_{22}\end{bmatrix} \begin{bmatrix} \underline{H}_{-j}&U_2 \\ 0 & D\end{bmatrix} 
+ \begin{bmatrix} \underline{H}_{-j}^H&0 \\U_2^H  & D^H\end{bmatrix}
 \begin{bmatrix} Y_{11} & Y_{12}\\Y_{12}^H & Y_{22}\end{bmatrix} \\
&\qquad\qquad -\begin{bmatrix} Z_j^H \\ \hat Z^H \end{bmatrix}BB^H \begin{bmatrix} Z_j~~\hat Z\end{bmatrix}
- \begin{bmatrix} h_j^H \\ U_1^H \end{bmatrix}\begin{bmatrix} h_j~~U_1\end{bmatrix}\\
 &= \begin{bmatrix} Y_{11}\underline{H}_{-j}  & Y_{11}U_2 + Y_{12}D \\Y_{12}^H\underline{H}_{-j} &Y_{12}^HU_2  +Y_{22}D\end{bmatrix} 
+ \begin{bmatrix} \underline{H}_{-j}^HY_{11}&  \underline{H}_{-j}^HY_{12}\\U_2^HY_{11}+D^HY_{12}^H  & U_2^HY_{12}+D^HY_{22}\end{bmatrix} \\
&\qquad\qquad -\begin{bmatrix} Z_j^HBB^HZ_j & Z_j^HBB^H\hat Z \\ \hat Z^HBB^HZ_j & \hat Z^HBB^H\hat Z \end{bmatrix}
- \begin{bmatrix} h_j^Hh_j & h_j^HU_1 \\ U_1^Hh_j & U_1^HU_1 \end{bmatrix}.
\end{align*}
The upper left block yields the Lyapunov equation for the unknown $Y_{11}$
\[
Y_{11} \underline{H}_{-j} + \underline{H}_{-j}^H Y_{11} - Z_{j}^HBB^HZ_{j}-h_{j}^Hh_{j} = 0.
\]
This is just \eqref{lyap_here}. The unique solution is given by $Y_j^{-1}$. Thus $Y_{11} = Y_{j}^{-1}.$ 

The upper right block yields the Sylvester equation for the unknown $Y_{12}$
\begin{equation}\label{eq:sylv}
0 = Y_j^{-1}U_2 + Y_{12}D +\underline{H}_{-j}^HY_{12} - Z_j^HBB^H\hat Z - h_j^HU_1.
\end{equation}
Owing to $D = \mu_{j+1}I_p$ we have $Y_{12}D = \mu_{j+1}Y_{12}.$ Hence, the solution $Y_{12}$ can be computed by solving a linear system
\begin{equation}\label{eq_y12}
(\mu_{j+1}I_{jp}+\underline{H}_{-j}^H)Y_{12} = -Y_j^{-1}U_2 + Z_j^HBB^H\hat Z + h_j^HU_1.
\end{equation}
Due to the assumption that all shifts lie in the right half plane, $\mu_{j+1}I_{jp}+\underline{H}_{-j}^H$ is nonsingular. Thus, the solution $Y_{12}$ is uniquely determined.

The lower right block yields the Lyapunov equation for the unknown $Y_{22}$
\begin{equation}\label{eq:lyap}
0 = Y_{12}^HU_2 + Y_{22}D + U_2^HY_{12} + D^HY_{22} - \hat Z^HBB^H\hat Z - U_1^HU_1.
\end{equation}
From this, $Y_{22}$ can be read off,
\begin{equation}\label{eq_y22}
Y_{22} = \frac{1}{2\text{Re}(\mu)}\left( -Y_{12}^HU_2 - U_2^HY_{12} + \hat Z^HBB^H\hat Z + U_1^HU_1\right).
\end{equation}

In summary, in order to determine $Y_{j+1}^{-1},$ it suffices to solve the linear system \eqref{eq_y12} and to compute $Y_{22}$ from \eqref{eq_y22}. In order to do so, we need to know $U_1$ and $U_2$ from \eqref{eq:brad} which depend on the choice of $\hat Z.$ 

It is possible to extend the BRAD by more than one block $\hat Z$ at a time. Assume that $\hat Z$ consists of $\ell$ blocks $\hat Z_j$, $\hat Z = \left[ \hat Z_1 ~~\hat Z_2 ~~\cdots~~\hat Z_\ell\right] \in \mathbb{C}^{n \times \ell p}.$ Then from \eqref{brad_here} with $Z_{j+\ell}$ instead of $Z_{j+1}$ (and appropriately adapted indices in the rest of the equation)  we have that \eqref{eq:brad} still holds with an upper triangular matrix $D \in \mathbb{C}^{p \times \ell p}, U_1 \in \mathbb{C}^{p\times \ell p}$ and $U_2 \in \mathbb{C}^{jp \times \ell p}.$ All derivations above still holds, such that $Y_{j+\ell}^{-1}$ can be computed as $\left[\begin{smallmatrix} Y_j^{-1} & Y_{12}\\ Y_{12}^H & Y_{22}\end{smallmatrix}\right]$ where $Y_{12}$ solves \eqref{eq:sylv}, while $Y_{22}$ solve \eqref{eq:lyap}. As before, in order to do so, we need to know $U_1$ and $U_2$ (as well as $D$) from \eqref{eq:brad} which depend on the choice of $\hat Z.$

In the next two subsections, we will consider two different possibilities for the choice of $\hat Z.$ Both choices allow for adding more than one shift at a time. The linear systems solves needed to supplement the block Krylov subspace accordingly can be solved simultaneously. This can be used in an efficient implementation to speed up the computations. The first choice allows for a faster computation of the approximations $X_j$ than any of the three equivalent algorithms from \cite{BenBKS18} in case $m$ is significantly larger than $p$  (and the Sherman-Morrision-Woodbury formula as well as a direct solver is used to solve the linear systems that occur), see Section \ref{sec_numex}.
The second choice for $\hat Z$ discussed  rediscovers the RADI algorithm \cite{BenBKS18} (up to some scaling). Hence, our approach gives a new interpretation of the RADI algorithm in terms of a BRAD and extends the known algorithm by allowing to add more than one shift at a time.

We conclude this section with a final remark which may be helpful for an efficient implementation of the algorithm.
\begin{remark}
Assume that $Y_j^{-1} = G_j^HG_j$ holds. Then $Y_{j+1}^{-1} = G_{j+1}^HG_{j+1}$ with
\[
G_{j+1} = \begin{bmatrix} G_j & G_j^{-H}Y_{12}\\ 0 & G_{22}\end{bmatrix}
\]
where $G_{22}$ is the Cholesky factor of
\[
Y_{22}-Y_{12}^HY_jY_{12} = G_{22}^HG_{22}.
\]
Hence, instead of $X_j = Z_jY_jZ_j^H$ we can consider $X_j = (Z_jG_j^{-1})(Z_jG_j^{-1})^H.$
\end{remark}
\subsection{Expanding the Krylov subspace by one or several blocks of the form $(A^H-\mu I_n)^{-1}R_j$}\label{subsec1}
A straightforward choice for $\hat Z$ is $\hat Z = (A^H-\mu I)^{-1}V_{j}T = (A^H-\mu I)^{-1}[C^H~~Z_j]T$ for some $T$ and $\mu \in \mathbb{C}_+\backslash \Lambda(A).$  Choosing $T = \left[\begin{smallmatrix} I_p\\ Y_jh_j^H\end{smallmatrix}\right]$ yields $\hat Z = (A^H-\mu I)^{-1}R_j,$ or, equivalently,
\begin{equation}\label{eq_linsys}
A^H\hat Z = R_j + \mu \hat Z = \begin{bmatrix}C^H~~Z_j~~\hat Z\end{bmatrix}
\begin{bmatrix} I_p\\Y_jh_j^H\\\mu I_p\end{bmatrix},
\end{equation}
as $R_j = C^H+  Z_jY_jh_j^H$ \eqref{eq:Rj}.
Thus, in the BRAD \eqref{eq:brad} we have $D = \mu I_p,$ $U_1=I_p$ and $U_2 = Y_jh_j^H$. This implies that \eqref{eq_y12} reduces to solving
\[
 (\mu I +\underline{H}_{-j}^H)Y_{12} = Z_j^HBB^H\hat Z,
\]
while \eqref{eq_y22} now reads
\[
Y_{22} = \frac{1}{2\text{Re}(\mu)}\left( -Y_{12}^H(Y_jh_j^H )- (Y_jh_j^H)^HY_{12} + \hat Z^HBB^H\hat Z + I_p\right).
\]
 Please note that only $Y_j^{-1}$ is available. Thus, in order to compute $Y_jh_j^H$ a linear system of equations with multiple right hand sides has to be solved.
The resulting algorithm is summarized in Algorithm \ref{alg2}. For ease of description, it is assumed that $j$ shifts are given and used in the given order. 

\begin{algorithm}[!ht]
\caption{Expanding the Krylov subspace by blocks of the form $\hat Z = (A^H-\mu I_n)^{-1}R_j$}\label{alg2}
\begin{algorithmic}[1]
\Require $A \in \mathbb{C}^{n \times n}, B\in \mathbb{C}^{n \times m}, C\in \mathbb{C}^{p \times n},$ and a set of shifts $\mathcal{S}_j = \{s_1, \ldots, s_j\}, s_i \in \mathbb{C}_+ \backslash \Lambda(A)$
\Ensure approximate solution $X_j = Z_jY_jZ_j^H$ of \eqref{ric}, residual factor $R_j$ such that $\mathcal{R}(X_j)=R_jR_j^H.$
\State Initialize $R_0 = C^H,h_1= I_p, $ and  $\underline{H}_{-1}=s_1 I_p.$
\State Solve $(A^H-s_1I_n) Z_1 = R_0.$ \Comment{$Z_1 = \hat Z.$}
\State Set $Y_1^{-1} = \frac{1}{2\text{Re}(s_1)}\left(Z_1^HBB^HZ_1 +I_p\right).$
\State Set $R_1 = C^H+Z_1Y_1h_1^H.$
\State Set $T_1 = B^HZ_1.$
\For {k = 2:j}
\State Solve $(A^H-s_kI_n)\hat Z = R_{k-1}.$
\State Solve  $(s_k I +\underline{H}_{-{k-1}}^H)Y_{12} = T_{k-1}^H(B^H\hat Z).$
\State Set $Y_{22} = \frac{1}{2\text{Re}(s_k)}\left( -Y_{12}^H(Y_{k-1}h_{k-1}^H )- (Y_{k-1}h_{k-1}^H)^HY_{12} + \hat Z^HBB^H\hat Z + I_p\right).$
\State Set $Y_k^{-1} = \left[\begin{smallmatrix}Y_{k-1}^{-1} &Y_{12}\\Y_{12}^H&Y_{22}\end{smallmatrix}\right].$
\State Set $Z_k = [ Z_{k-1}~~\hat Z].$
\State Set $h_k=[h_{k-1}~~I_p].$
\State Set $R_k = C^H+Z_kY_kh_k^H.$
\State Set $\underline{H}_{-k}=\left[\begin{smallmatrix} \underline{H}_{-(k-1)}&Y_{k-1}h_{k-1}^H\\ 0 & s_k I_p\end{smallmatrix}\right]$
\State Set $T_k = [T_{k-1}~~B^H\hat Z].$
\EndFor
\end{algorithmic}
\end{algorithm}

Next we consider the case of extending the BRAD by several blocks of the form $(A^H-\mu I)^{-1}R_j$ at once.
Assume that  shifts $\mu_1, \ldots, \mu_\ell  \in \mathbb{C}_+ \backslash \Lambda(A)$ are given. Then the linear systems
\[
\hat Z_i = (A^H-\mu_iI_n)^{-1}R_j, \qquad i = 1, \ldots, \ell
\]
are independent of each other and can be solved at the same time (in parallel). Each of these linear systems can be written in the form \eqref{eq_linsys}. Thus the expanded BRAD \eqref{eq:brad} with $\hat Z = \left[ \hat Z_1~~\hat Z_2~~\cdots ~\hat Z_\ell\right]$ is given by
\begin{align*}
A^H\left[ C^H~~Z_j~~\hat Z\right] \begin{bmatrix}0\\I_{(j+\ell)p}\end{bmatrix} &=
\left[ C^H~~Z_j~~\hat Z\right] 
{\footnotesize \begin{bmatrix}
h_j & I_p &I_p & \cdots & I_p\\
\underline{H}_{-j} & Y_jh_j^H & Y_jh_j^H & \cdots & Y_jh_j^H\\
0 & \mu_1 I_p & 0 & \cdots & 0\\
0 & 0 & \mu_2 I_p &  \ddots & \vdots\\
\vdots& \vdots &\ddots &\ddots &0\\
0 &0&\cdots &0& \mu_\ell I_p
\end{bmatrix}}\\
& = \left[ C^H~~Z_j~~\hat Z\right] 
 \begin{bmatrix}
h_j & U_1\\
\underline{H}_{-j} & U_2\\
0 &D
\end{bmatrix} 
\end{align*}
with $U_1 = [I_p~~I_p~~\cdots~I_p]$, $D= \operatorname{diag}(\mu_1I_p, \ldots, \mu_\ell I_p)$ and $U_2 =Y_jh_j^HU_1.$ From the definition \eqref{brad} of a BRAD it follows that the BRAD has been expanded with the shifts $\mu_1, \ldots, \mu_\ell.$ Following the derivation of the algorithm in the previous section, Algorithm \ref{alg2} has to be modified slightly in order to take this modification into account
\begin{itemize}
\item solve $(A^H-\mu_iI_n) \hat Z_i = R_j$ for $i = 1, \ldots, \ell,$
\item set $\hat Z = \left[ \hat Z_1~~\hat Z_2~~\cdots ~\hat Z_\ell\right],$ $U_1 = [I_p~~I_p~~\cdots~I_p]$, $D= \operatorname{diag}(\mu_1I_p, \ldots, \mu_\ell I_p)$ and $U_2 =Y_jh_j^HU_1,$
\item solve \eqref{eq:sylv} for $Y_{12} \in \mathbb{C}^{jp \times \ell p}$  and \eqref{eq:lyap} for $Y_{22} \in \mathbb{C}^{\ell p \times \ell p},$
\item set $Y_k^{-1} = \left[\begin{smallmatrix}Y_{k-1}^{-1} &Y_{12}\\Y_{12}^H&Y_{22}\end{smallmatrix}\right]$
and $Z_k = [ Z_{k-1}~~\hat Z],$
\item set $h_k=[h_{k-1}~~U_1]$ and $\underline{H}_{-k}=\left[\begin{smallmatrix} \underline{H}_{-(k-1)}&Y_kh_k^HU_1\\ 0 & D\end{smallmatrix}\right],$
\item set $R_k = C^H+Z_kY_kh_k^H.$
\end{itemize}
Please note that the first step can be performed in parallel which may allow for an efficient fast  implementation. 

\begin{remark}\label{rem_review2}
In Algorithm 2 in \cite{LinS15} as well as in the equivalent Algorithm 1 in \cite{MasOR16} the Krylov subspace is advanced in each iteration step by a block involving $(A^H-\mu_iI_n)^{-1}$ similar to our approach considered here. For the following discussion in this remark we will make use of the notation in \cite{MasOR16}.  The first block $V_1$ used is essentially  the same one as in our Algorithm \ref{alg2} ($V_1 = -Z_1$). In particular, Line 1 of Algorithm 1 in \cite{MasOR16} gives $A^HV_1 = \alpha_1 V_1 - C^H.$  The following blocks are of the form $V_j =V_{j-1} -(\alpha_{j}+\overline{\alpha_{j-1}})(\alpha_jI-A^H)^{-1}V_{j-1}.$ 
Some rewriting of this expression yields  $A^HV_j = \alpha_j V_j + A^H V_{j-1} + \overline{\alpha_{j-1}}V_{j-1}.$
Consider this for $j=2$ and insert the expression for $A^HV_1.$ This gives
\[
A^HV_2 = \alpha_2V_2 +\operatorname{Re}(\alpha_1)V_1 -C^H
= [C^H~~V_1~~V_2]\begin{bmatrix} -I\\ \operatorname{Re}(\alpha_1)I\\\alpha_2I\end{bmatrix}.
\]
In this fashion, we were able to find  $U_2$ such that
\[
A^HV_j = \left[C^H~~[V_1~\cdots~V_{j-1}] ~~V_j\right]
\begin{bmatrix} -I\\U_2\\\alpha_j I\end{bmatrix} 
\]
holds. Thus, Algorithm 1 from \cite{MasOR16} and Algorithm 2 from  \cite{LinS15}  do fit into our BRAD-approach as \eqref{eq_linsys} holds.
\end{remark}

\subsection{Expanding the Krylov subspace by one or several blocks of the form $(A^H-X_jBB^H -\mu I_n)^{-1}R_j$}\label{subsec2}
Another possible choice for $\hat Z$ is $(A^H-X_jBB^H-\mu I)^{-1}R_j$.  This gives
\begin{equation}\label{eq_linsys2}
A^H\hat Z = R_j + X_jBB^H\hat Z+ \mu \hat Z = \begin{bmatrix}C^H~~Z_j~~\hat Z\end{bmatrix}
\begin{bmatrix} I_p\\Y_jh_j^H+ Y_jZ_j^HBB^H\hat Z\\\mu I_p\end{bmatrix}
\end{equation}
as $R_j =\left[C^H~~Z_j\right]\left[\begin{smallmatrix}I_p\\Y_jh_j^H\end{smallmatrix}\right]$ \eqref{eq:Rj} and $X_j = Z_jY_jZ_j^H.$
Thus, in the BRAD \eqref{eq:brad} we have $U_1=I_p$ and $U_2 = Y_jh_j^HU_1+Y_jZ_j^HBB^H\hat Z.$ This implies that \eqref{eq_y12} reduces to
\[
 (\mu I +\underline{H}_{-j}^H)Y_{12} =0,
\]
which gives $Y_{12}=0.$
The Lyapunov equation \eqref{eq_y22} simplifies to
\begin{equation}\label{y22_hier}
Y_{22} = \frac{1}{2\text{Re}(\mu)}\left( I_p+ \hat Z^HBB^H\hat Z \right).
\end{equation}
The residual factor  $R_ {j+1} = C^H + Z_{j+1}Y_{j+1}h_{j+1}^H$ is updated as follows
\[ R_{j+1} = \left[C^H~~Z_{j+1}\right] \begin{bmatrix} I_p\\Y_{j+1}h_{j+1}^H\end{bmatrix} = \left[C^H~~Z_{j}~~\hat Z\right] \begin{bmatrix} I_p\\Y_{j}h_{j}^H\\Y_{22}^{-1}U_1^H\end{bmatrix} = R_j + \hat Z Y_{22}^{-1}U_1^H,
\]
as $Y_{j+1}^{-1} = \left[\begin{smallmatrix}Y_j^{-1}&0\\0& Y_{22}\end{smallmatrix}\right].$
 Please note that only $Y_j^{-1}$ is available. Thus, in order to compute $Y_jh_j^H$ a linear system of equations with multiple right hand sides has to be solved.
The resulting algorithm is summarized in Algorithm \ref{alg3}. As before, for ease of description, it is assumed that $j$ shifts are given and used in the given order. 

\begin{algorithm}
\caption{Expanding the Krylov subspace by blocks of the form $\hat Z = (A^H-X_jBB^H-\mu I_n)^{-1}R_j$}\label{alg3}
\begin{algorithmic}[1]
\Require $A \in \mathbb{C}^{n \times n}, B\in \mathbb{C}^{n \times m}, C\in \mathbb{C}^{p \times n},$ and a set of shifts $\mathcal{S}_j = \{s_1, \ldots, s_j\}, s_i \in \mathbb{C}_+ \backslash \Lambda(A)$
\Ensure approximate solution $X_j = Z_jY_jZ_j^H$ of \eqref{ric}, residual factor $R_j$ such that $\mathcal{R}(X_j)=R_jR_j^H.$
\State Initialize $R_0 = C^H, K_0 = 0 ,Y_0^{-1} =[]$ and  $Z_0=[].$
\For {k = 1:j}
\State Solve $(A^H-K_{k-1}B^H-s_kI_n)\hat Z = R_{k-1}.$
\State Set $Y_{22} = \frac{1}{2\text{Re}(s_k)}\left( \hat Z^HBB^H\hat Z + I_p\right).$
\State Set $Y_k^{-1} = \left[\begin{smallmatrix}Y_{k-1}^{-1} &\\&Y_{22}\end{smallmatrix}\right].$
\State Set $Z_k = [ Z_{k-1}~~\hat Z].$
\State Set $R_k = R_{k-1}+ \hat Z Y_{22}^{-1}.$
\State Set $K_k = K_{k-1}+\hat Z Y_{22}^{-1}(\hat Z^HB).$
\EndFor
\end{algorithmic}
\end{algorithm}

\begin{remark}\label{rem2_review2}
Although the approach taken in \cite{BenBKS18} to derive an algorithm for solving the algebraic Riccati equation \eqref{ric0} is
quite different from that taken here, Algorithm \ref{alg3} is essentially the same as  Algorithm 1 in \cite{BenBKS18}. They differ only by the scaling of the matrices involved.

 The algorithm as derived in \cite{BenBKS18} allows for a nonzero initial guess $X_0$. This may be appropriate in case $A$ is not stable. In that case, the equivalence \eqref{eq:equiv1} does not hold. Our approach does not allow for a nonzero initial guess, as $X_0 =0$ and $R_0 =C^H$ is build into Algorithm \ref{alg1} - \ref{alg3} by the choice of the first block $C^H$ in the underlying Krylov subspace. But, as noted in \cite[Theorem 1(b)]{BenBKS18}, if $\tilde X$ is the solution of
\begin{equation}\label{eq1aa}
\tilde A^H \tilde X +\tilde X\tilde A + \tilde Q - \tilde X BB^H \tilde X = 0
\end{equation}
with $\tilde A = A-BB^HX_0$ and $\tilde Q = \mathcal{R}(X_0),$ then $X = X_0+\tilde X$ is a solution to the original CARE \eqref{ric0}. Thus, solving \eqref{eq1aa} with a suitable $X_0$ allows for using a nonzero initial starting guess/stabilizing $A.$

As pointed out in \cite{BenBKS18,BenBKS20}, if the CARE  \eqref{ric0} arises in LQR-optimal control, then only the feedback gain $K=XB$ is needed. RADI and hence also Algorithm \ref{alg3} can operate on approximate gains $K_j$ alone without storing the whole low-rank factor $Z_j$. This is due to the fact that $Y_j$ is a block diagonal matrix. This is not possible in Algorithm \ref{alg2} in which $Y_j$ has no special form which can be exploited to obtain a direct update of $K_j$ using only recent information.
\end{remark}

\begin{remark}\label{rem:SMW}
In general, the matrix $A^H-K_jB^H-\mu I$ with $K_j = X_jB$ is a dense matrix even if $A$ is sparse. However, as $K_jB^H$ is of rank $m$, it is proposed in \cite[Section 4,2]{BenBKS18} to use the Sherman-Morrison-Woodbury formula to speed up computations. With $L = (A^H-\mu I)^{-1}R_j$ and $N =  (A^H-\mu I)^{-1}K_j$ the formula reads
\[
(A^H-K_jB^H-\mu I)^{-1}R_j = L+N(I_m-B^HN)^{-1}B^HL.
\]
Thus, in case $A$ is sparse, first one large-scale sparse linear system with $p+m$ right-hand sides has to be solved in order to determine $L$ and $N$, then a small $m \times m$ possibly dense system has to be solved. This approach allows for solving the resulting  linear systems by a direct solver.
\end{remark}

\begin{remark}\label{rem_lyap_adi}
As already mentioned in Remark \ref{rem_lyap}, Algorithms \ref{alg1} - \ref{alg3} can be applied to continuous Lyapunov equations \eqref{eq_lyapsolv}. In this case, Algorithm \ref{alg2} and Algorithm \ref{alg3} become identical and simplify considerably boiling down to the well-known low-rank ADI iteration for Lyapunov equations (see, e.g., \cite{BenKS13a,BenKS13,LiW02,Pen00} and the reference therein)
\begin{itemize}
\item solve $(A^H-s_kI_n)\hat Z = R_{k-1},$
\item set $Z_k = [ Z_{k-1}~~\sqrt{2\text{Re}(s_k)}\hat Z],$
\item set $R_k = R_{k-1}+ 2\text{Re}(s_k)\hat Z,$
\end{itemize}
with $X_j = Z_jZ_j^H.$
\end{remark}

To conclude the discussion, we consider the case of extending the BRAD by several blocks of the form $(A^H-K_jB^H-\mu I)^{-1}R_j$ at once.
Assume that shifts $\mu_1, \ldots, \mu_\ell  \in \mathbb{C}_+ \backslash \Lambda(A)$ are given. Then the linear systems
\[
\hat Z_i = (A^H-K_jB^H-\mu_iI_n)^{-1}R_j, \qquad i = 1, \ldots, \ell
\]
are independent of each other and can be solved at the same time (in parallel). Each of these linear systems can be written in the form \eqref{eq_linsys2}. Thus the expanded BRAD \eqref{eq:brad} with $\hat Z = \left[ \hat Z_1~~\hat Z_2~~\cdots ~\hat Z_\ell\right]$ is given by
\[
A^H\left[ C^H~~Z_j~~\hat Z\right] \begin{bmatrix}0\\I_{(j+\ell)p}\end{bmatrix} =
 \left[ C^H~~Z_j~~\hat Z\right] 
 \begin{bmatrix}
h_j & U_1\\
\underline{H}_{-j} & U_2\\
0 &D
\end{bmatrix} 
\]
with $U_1 = [I_p~~I_p~~\cdots~I_p]$, $D= \operatorname{diag}(\mu_1I_p, \ldots, \mu_\ell I_p)$ and $U_2 =Y_jh_j^HU_1+Y_jZ_j^HBB^H\hat Z.$ As before from the definition \eqref{brad} of a BRAD it follows that the BRAD has been expanded with the shifts $\mu_1, \ldots, \mu_\ell.$ Algorithm \ref{alg3} has to be modified slightly in order to take this modification into account
\begin{itemize}
\item solve $(A^H-K_jB^H-\mu_iI_n) \hat Z_i = R_j$ for $i = 1, \ldots, \ell,$
\item set $\hat Z = \left[ \hat Z_1~~\hat Z_2~~\cdots ~\hat Z_\ell\right],$ $U_1 = [I_p~~I_p~~\cdots~I_p]$, $D= \operatorname{diag}(\mu_1I_p, \ldots, \mu_\ell I_p)$ and $U_2 =(Y_jh_j^H+Y_jZ_j^HBB^H\hat Z)U_1,$
\item solve \eqref{eq:lyap} for $Y_{22} \in \mathbb{C}^{\ell p \times \ell p}$ (making use of $Y_{12}=0$),
\item set $Y_k^{-1} = \left[\begin{smallmatrix}Y_{k-1}^{-1} &0\\0&Y_{22}\end{smallmatrix}\right]$
and $Z_k = [ Z_{k-1}~~\hat Z],$
\item set $R_k = R_{k-1}+\hat Z Y_{22}^{-1}U_1^H,$
\item set $K_k = K_{k-1}+\hat ZY_{22}^{-1}(\hat Z^HB).$
\end{itemize}
Please note that the first step can be performed in parallel which may allow for an efficient fast  implementation. This extension of Algorithm \ref{alg3}/the RADI algorithm \cite[Algorithm 1]{BenBKS18} is new.  The parallelization follows easily from the setting considered here, but is not obvious from the context discussed in \cite{BenBKS18}.

\begin{remark}
The parallelization can also be incorporated into the low-rank ADI iteration for Lyapunov equations which is just Algorithm \ref{alg3}, see Remark \ref{rem_lyap_adi}.
See \cite[Remark 5.23]{Wol14} for a different approach for the parallelization of the ADI iteration for Lyapunov equations.
\end{remark}

\subsection{Realification in case of real matrices $A,B,C$}\label{subsec:real}
In case of real system matrices $A, B, C$ and two complex conjugate shifts a modification of the algorithm makes sure that the iterates remain real valued. The idea presented below is based on the double vector variant described in \cite[Section 3]{Ruh10} as well as on \cite[Section 4.3]{BenBKS18}. The use of complex arithmetic can not be completely avoided. 
We will discuss the procedure for the choice $\hat Z = (A^T-\mu I_n)^{-1}R_j.$ It can be adapted for the choice $\hat Z = (A^T-K_jB^H -\mu I_n)^{-1}R_j$ in a straightforward way.

Let $R_j$ have only real entries. Let the BRAD \eqref{brad_here} be expanded by the two conjugate complex shifts $\mu, \overline{\mu} \in \mathbb{C} \backslash \Lambda(A)$, $\mu = a+\imath b, a, b \in \mathbb{R}, b \neq 0.$ Let $W = (A^T-\mu I_n)^{-1}R_j.$ Then $\overline{W} =  (A^T-\overline{\mu} I_n)^{-1}R_j.$ With
\[
S = \frac{1}{2}\begin{bmatrix}I_p & -\imath I_p\\ I_p & \imath I_p\end{bmatrix} 
\]
we have
\[
\left[ W~~\overline{W}\right] S = \left[ \operatorname{Re}(W)~~ \operatorname{Im}(W)\right] \qquad \text{and}\qquad 
S^{-1}\begin{bmatrix} \mu I_p\\&\overline{\mu}I_p\end{bmatrix} S = \begin{bmatrix} aI_p & b I_p\\ -b I_p& a I_p\end{bmatrix}.
\]
Expanding the BRAD  \eqref{eq:brad} with the complex blocks $W$ and $\overline{W}$ yields
\[
A^T\left[ C^T~~Z_j~~W~~\overline{W}\right] \begin{bmatrix}0\\I_{(j+2)p}\end{bmatrix} =
\left[ C^T~~Z_j~~W~~\overline{W}\right] 
 \begin{bmatrix}
h_j & I_p &I_p \\
\underline{H}_{-j} & Y_jh_j^H & Y_jh_j^H \\
0 & \mu I_p & 0 \\
0 & 0 & \overline{\mu} I_p 
\end{bmatrix}.
\]
Transforming this BRAD with the matrix $\left[\begin{smallmatrix} I_{(j+1)p} & 0\\ 0 & S\end{smallmatrix}\right]$ gives
\begin{align}
\begin{split}
&A^T\left[ C^T~~Z_j~~\operatorname{Re}(W)~~\operatorname{Im}(W)\right] \begin{bmatrix}0\\I_{(j+2)p}\end{bmatrix} \\
&\qquad =
\left[ C^T~~Z_j~~\operatorname{Re}(W)~~\operatorname{Im}(W)\right] 
 \begin{bmatrix}
h_j & I_p &0\\
\underline{H}_{-j} & Y_jh_j^H & 0 \\
0 & a I_p & b I_p \\
0 & -b I_p & a I_p 
\end{bmatrix}.
\end{split}\label{nobrad}
\end{align}
The Krylov basis is expanded with the real block vectors  $\operatorname{Re}(W)$ and $\operatorname{Im}(W).$
Note that only one (complex) system solve is necessary for the expansion with the complex conjugate pair of
shifts $\mu$ and $\overline{\mu}.$ 

The relation \eqref{nobrad} is no longer a BRAD as the rightmost matrix is no longer a block upper Hessenberg matrix. This can be fixed by transforming \eqref{nobrad} with the matrix $\left[\begin{smallmatrix} I_{(j+1)p} & 0\\ 0 & P\end{smallmatrix}\right]$ for \[P = [e_1~e_{p+1}~e_2~e_{p+2}~\cdots~e_p~e_{2p}].\]
This corresponds to permuting the columns of $\operatorname{Re}(W)$ and $\operatorname{Im}(W)$ so that the Krylov basis is expanded by $\left[ \operatorname{Re}(w_1)~\operatorname{Im}(w_1)~
\cdots \operatorname{Re}(w_p)~\operatorname{Im}(w_p)\right]$ for $W = \left[w_1~w_2~\cdots~w_p\right].$ It transforms the
 lower right block $
\left[\begin{smallmatrix} a & b \\ -b & a \end{smallmatrix}\right] \otimes I_p$ into $I_p \otimes \left[\begin{smallmatrix} a & b \\ -b & a \end{smallmatrix}\right]$ which is  a block diagonal matrix with $2 \times 2$ blocks on the diagonal,
\begin{align*}
&A^T\left[ C^T~~Z_j~~\operatorname{Re}(w_1)~\operatorname{Im}(w_1)~
\cdots \operatorname{Re}(w_p)~\operatorname{Im}(w_p)\right] \begin{bmatrix}0\\I_{(j+2)p}\end{bmatrix} =\\
&
\hspace*{1cm}\left[ C^T~~Z_j~~\operatorname{Re}(w_1)~\operatorname{Im}(w_1)~
\cdots \operatorname{Re}(w_p)~\operatorname{Im}(w_p)\right] 
\left[ \begin{array}{c|cc}
h_j & I_p &I_p \\
\underline{H}_{-j} & Y_jh_j^H & Y_jh_j^H \\ \hline
0 & \multicolumn{2}{c}{I_p \otimes \left[\begin{smallmatrix} a & b \\ -b & a \end{smallmatrix}\right]}
\end{array}\right]
\end{align*}
which is a BRAD again.
In case only complex conjugate pairs of shifts are used in the fashion discussed above, the right most matrix will be a quasi upper triangular matrix which can be used to solve the Sylvester equation \eqref{eq:sylv} efficiently with established software packages.

\subsection{Choice of Shifts}
The methods proposed in this work require shift parameters to achieve a rapid convergence just like the four equivalent Riccati ADI algorithms. As the proposed method are equivalent to the RADI algorithm, we point the readers to, e.g.,
\cite{BenBKS18,BenBKS20,Kue16,Sim16} for different suggestions on how to choose the shifts and comparisons of different approaches. 
The interplay of the possible parallelization and the choice of the shifts concerning the convergence behavior may be crucial, but an in-depth discussion of this aspect is beyond this work.

{\subsection{Multiple Use of the Same Shift}\label{subsec_review2}
In Algorithm \ref{alg2} and Algorithm \ref{alg3} (as well as in extending the BRAD) by several blocks at once the same shift can be used multiple times. Although, compared to distinct shifts in each step, this may hinder convergence, a measurable savings in calculation time might be achieved in case a (sparse) direct solver is used in order to solve the linear systems, see, e.g., \cite[Section 4]{Kue19} for a discussion on this aspect in the context of low-rank ADI solvers for Lyapunov equations.
In case a shift is used more than once, \eqref{eq:same} no longer applies, the definition of $\mathcal{L}_j$ needs to be adapted appropriately.

\section{Numerical Experiments}\label{sec_numex}
In this section we compare Algorithms \ref{alg2} and \ref{alg3}. Please note that Algorithm \ref{alg3} is just the RADI algorithm from \cite{BenBKS18}, the only difference is on how the necessary scaling (in terms of $2\text{Re}(s_k)$) is incorporated. 
Algorithms \ref{alg2} and \ref{alg3} compute exactly the same iterates $X_j$ (when using the same set of shifts and $X_0=0$), as we proved in Theorem \ref{theo:main}. Moreover, another consequence of Theorem \ref{theo:main} is that $X_j = X_j^\text{qadi} = X_j^\text{cay} = X_j^\text{radi}$ holds. An extensive comparison of the low-rank qADI algorithm, the
Cayley transformed subspace iteration, and the RADI iteration has been presented in \cite{BenBKS18,BenBKS20}. 
We complement those findings by comparing Algorithms \ref{alg2} and \ref{alg3} with respect to their timing performance. 

The most time consuming part of both algorithms is solving the $n \times n$ system of linear equations at the beginning of each iteration step. While the linear system in Algorithm \ref{alg2} is sparse, the one in Algorithm \ref{alg3} is in general dense. To make our implementation comparable with the RADI implementation from \cite{BenBKS18}, our implementation makes use of  the Sherman-Morrison-Woodbury formula as discussed in Remark \ref{rem:SMW} in order to rewrite the dense linear systems in Algorithm \ref{alg3} into sparse linear systems. Thus in each iteration step of Algorithm \ref{alg3} one sparse large-scale system with $m + p$ right-hand sides has to be solved, while in Algorithm \ref{alg2} just one sparse large-scale system with $p$ right-hand sides has to be handled per iteration step.
Hence, it is to be expected that Algorithm \ref{alg2} may be faster than Algorithm \ref{alg3} in case $m$ is significantly larger than $p$.
To look at this aspect in more detail, we have chosen the following test structure: First, the shifts are calculated so that both algorithms perform the same number of iteration steps with the same shifts. This part has not been included in the timings reported. In each iteration step, the time required to solve the respective linear system of equations is measured. In addition, the total time needed by the algorithms for the required number of iterations steps is noted.

Algorithms \ref{alg2} and \ref{alg3} have been implemented in MATLAB including
 realification in case of a complex shift as explained in Section \ref{subsec:real} as well as the modification needed to handle  generalized Riccati equations \eqref{ric_gen} with an additional system matrix $E$ as discussed in Section \ref{subsec:genR}. 
All linear systems are solved by a direct solver via MATLAB's $\backslash$-operator. The shifts are precomputed using the efficient implementation of the RADI algorithm in the MATLAB toolbox M.E.S.S.-2.2 \cite{SaaKB22} (employing default settings for \texttt{mess\_lrradi} which implies that the shift strategy residual Hamiltonian shifts \cite[Section 4.5.1]{BenBKS18} is used). 
The experimental code used to generate the results presented can be found at \cite{zenodo}.
All experiments are performed in MATLAB R2023a on an Intel(R) Core(TM) i7-8565U CPU @ 1.80GHz
1.99 GHz with 16GB RAM. 

The first example considered is the well-known steel profile cooling model
from the Oberwolfach Model Reduction Benchmark Collection \cite{morwiki_steel,BenS05b}. This example (often termed RAIL)  comes in different problem sizes $n$, but fix $m =7$ and $p = 6$. We used the one with $n = 79,841.$ The system matrices $E$ and $A$ are symmetric positive and negative definite, resp.. All (with \texttt{mess\_lrradi}) precomputed shifts are real.
The plot on the left-hand side in Figure \ref{fig1} displays the computational times measured  for the linear system solve in each iteration step. Usually, the system solves in Algorithm \ref{alg2} need less time than those in Algorithm \ref{alg3}. In Table \ref{tab11}, the total computational time for solving the linear systems as well as the computational time for the entire algorithms is given. It can be seen that Algorithm \ref{alg2} is slightly faster than Algorithm \ref{alg3} in both of these aspects. The impact of the larger number of right-hand sides in the system solves in Algorithm \ref{alg3} compared to Algorithm \ref{alg2} comes only little to bear here as $m$ is fairly small.

The second example considered is the convection-diffusion benchmark example from MORwiki - Model Order Reduction Wiki \cite{morwiki,lyapack}. The examples are constructed with 

\begin{center}
\begin{minipage}{12cm}
\texttt{A = fdm\_2d\_matrix(100,'10*x','100*y','0');}\\
\texttt{B = fdm\_2d\_vector(100,'.1<x<=.3');}\\
\texttt{ C = fdm\_2d\_vector(100,'.7<x<=.9')'; }\\
\texttt{E = speye(size(A));} 
\end{minipage}
\end{center}
resulting in a SISO system of order $n = 10,000.$ Among the $46$ (with \texttt{mess\_lrradi}) precomputed shifts there are $20$ real ones and $13$ pairs of complex-conjugate shifts. The plot on the right-hand side in Figure \ref{fig1} displays the computational times measured  for the linear system solve in each iteration step. As can be seen, real shifts have been used in the iteration steps 1--3, 8, 10, 12--15, 18--19, 21--22, 25--26, 29, and 32. The steps associated with complex shifts are more expansive than those with real shifts. Overall, Algorithm \ref{alg2} and Algorithm \ref{alg3} perform alike in terms of computational time.

\begin{figure}[htb!]
\begin{center}
\includegraphics[scale=0.45]{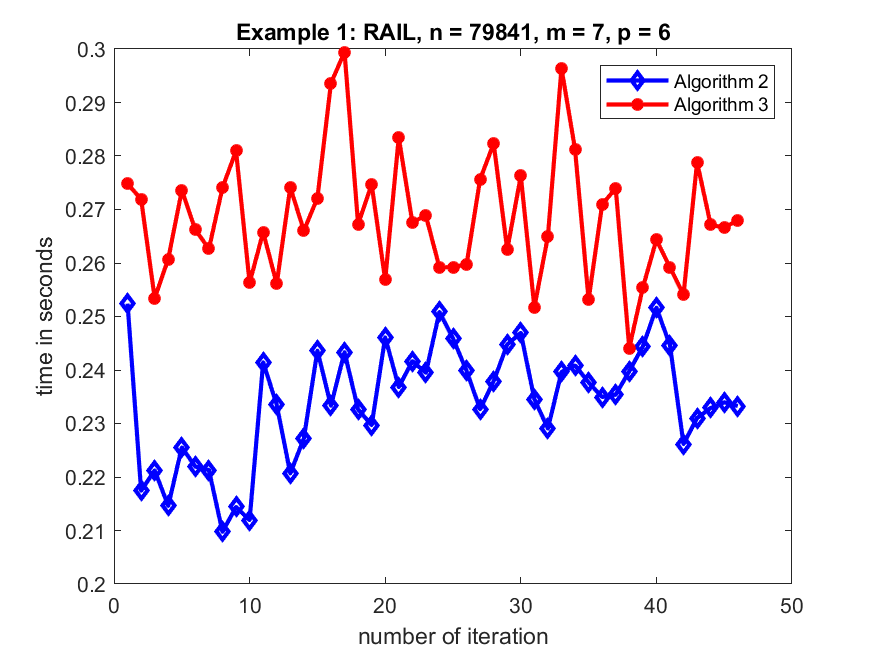}  \includegraphics[scale=0.45]{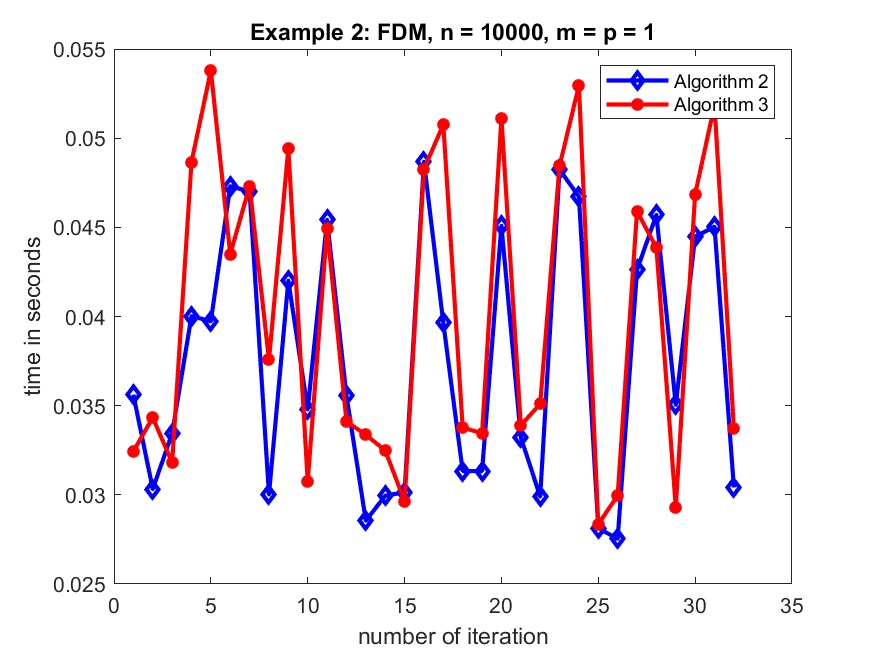}
\end{center}
\caption{Computational time for the linear system solve in each iteration step for Examples 1 and 2.\label{fig1}}
\end{figure}

\begin{table}[htb!]
\begin{center}
\begin{tabular}{l|c|rrr|rrr}
&& \multicolumn{3}{c|}{Algorithm 2}& \multicolumn{3}{c}{Algorithm 3}\\
&n& lin. solves & misc & total & lin. solves & misc & total \\ \hline
Ex. 1 & $79,841$ & $10.7688$ & $0.8682$ &   $11.6370$ &  $12.3457$ & $0.8005$ & $13.1462$\\
Ex. 2 & $10,000$ & $1.2030$ & $0.1053$ &   $1.3083$ &   $1.2817$ & $0.0308$ &  $1.3125$
\end{tabular}
\end{center}
\caption{Computational time in seconds for different parts of the algorithms for Examples 1 and 2.\label{tab11}}
\end{table}

While in the first two examples both $m$ and $p$ were either identical or differed only by one, in our third example we will consider $m$ significantly larger than $p$. In this case, Algorithm \ref{alg2} can show its potential for problems with many more inputs than outputs. We consider the matrix lung2 available from The SuiteSparse Matrix Collection\footnote{\url{https://sparse.tamu.edu}}  (formerly known as the University of Florida Sparse Matrix Collection) via the matrix ID 894 \cite{DavH11}, modelling processes in the human lung. 
We employ this example with the negated system matrix $-A \in \mathbb{R}^{109460 \times 109460}, E = I$ and randomly chosen $C^H \in \mathbb{R}^{109460 \times 3}, B \in \mathbb{R}^{109460\times m}, m = 3k$ (using \texttt{sprandn} with a density of 0.1). While for $p=m$ (that is, $k=1$) Algorithm \ref{alg3} has a faster overall run time than Algorithm \ref{alg2}, as soon as $k$ (and hence $m$) increases, Algorithm \ref{alg2} is faster than Algorithm \ref{alg3} as can be seen from the data given in Table \ref{tab2}. Recall, that while the number of right-hand sides for each sparse large-scale system solve is just $p$ for Algorithm \ref{alg2}, there are $m + p = kp+p =(k+1)p$ right-hand sides for each such system solve in Algorithm \ref{alg3}. The timings for Algorithm \ref{alg2} are more dependent on the number of shifts chosen than on $m$, while the timings for Algorithm \ref{alg3} depend on both. The larger $m$ is compared to $p$, the better Algorithm \ref{alg2} performs in terms of computational time.

%

\begin{table}[hbt!]
\begin{center}
\begin{tabular}{l|c|rrr|rrr}
& no. of& \multicolumn{3}{c|}{Algorithm 2}& \multicolumn{3}{c}{Algorithm 3}\\
p = 3& shifts& lin. solves & misc & total & lin. solves & misc & total \\ \hline
m = p&$106$  &$17.3222$ & $8.8966$ &   $26.2188$ &  $22.0616$ & $1.1428$ & $23.2044$\\
m = 5p  &$89$&  $13.6938$ & $6.0072$ &   $19.7010$ &   $29.1169$ & $1.2695$ &  $30.3864$\\
m = 30p &$80$ & $14.0496$ & $6.2356$ &   $20.2852$ &   $105.6359$ & $2.7577$ &  $108.3936$\\
m = 100p &$72$ & $12.2363$ & $5.2020$ &   $17.4383$ &   $385.9790$ & $6.0403$ &  $392.0193$
\end{tabular}
\end{center}
\caption{Computational time in seconds for different parts of the algorithms for Example 3 with varying $m$ and fixed $p$.\label{tab2}}
\end{table}

\section{Concluding Remarks}
In this paper, we have suggested a new family of low-rank algorithms for computing
solutions of large scale Riccati equations based on a block rational Arnoldi decomposition and an associated block rational Krylov subspace spanned by $A^H$ and $C^H.$ We have shown that these algorithms produce
exactly the same iterates as the RADI algorithm \cite{BenBKS18} (and three other previously known methods) (when using the same set of parameters). 
We have suggested two specific versions of the general algorithm; one turns out to be equivalent to the RADI algorithm, the other one yields a computationally more efficient way to generate the approximate solutions $X_j$ than the RADI algorithm as well as the other previously known equivalent methods  in case $m$ is significantly larger than $p$  (in case the Sherman-Morrision-Woodbury formula and a direct solver is used to solve the linear systems that occur).  In case the linear systems are solved by any other means this advantage might disappear.
The general approach allows for adding more than one shift at a time, so that a number of the linear systems to be solved can be solved simultaneously. A discussion of the possible parallelization when adding more than one shift at a time and the choice of shifts in such a case is beyond the scope of this paper. 

\section*{Acknowledgements}
The authors  thank the reviewers for the exceptionally careful reading of the first draft of this paper and the many critical and very helpful comments which helped us to significantly improve the presentation. In particular, Remark \ref{rem_review2}, most of Remark \ref{rem2_review2} and Section \ref{subsec_review2} are due to one of the reviewers.


\end{document}